\documentclass[11pt]{article}
\usepackage{graphicx}
\usepackage{textcomp}
\usepackage{amsmath,amsfonts,amscd,amsthm,amssymb}
\usepackage{epsfig}
\usepackage{makeidx}
\makeindex
\newcommand{\ncom}{\newcommand}
\newcommand{\be}{\begin{eqnarray*}}
\newcommand{\ee}{\end{eqnarray*}}
\newcommand{\ben}{\begin{eqnarray}}
\newcommand{\een}{\end{eqnarray}}
\newcommand{\dps}{\displaystyle}
\newcommand{\tl}{\tilde}

\newcommand{\pr}{\partial}
\newcommand{\nab}{\nabla}

\newcommand{\bu}{{\bf u}}
\newcommand{\bU}{{\bf U}}

\newcommand{\bw}{{\bf w}}

\newcommand{\e}{{\bf e}}

\newcommand{\bJ}{{\bf J}}

\newcommand{\bv}{{\bf v}}

\newcommand{\bta}{\mbox{\boldmath $\eta$}}
\newcommand{\bzeta}{\mbox{\boldmath $\zeta$}}

\newcommand{\bphi}{\mbox{\boldmath $\phi$}}
\newcommand{\brho}{\mbox{\boldmath $\rho$}}

\newcommand{\bxi}{\mbox{\boldmath $\xi$}}
\newcommand{\bH}{{\bf {H}}}

\newcommand{\bL}{\bf {L}}

\newcommand{\f}{{\bf {f}}}
\newcommand{\bvh}{{\bf v_h}}

\newtheorem{tdf}{Theorem}[section]

\newtheorem{ldf}{Lemma}[section]
 
\ncom{\ul}{\underline}
\ncom{\beq}{\begin{equation}}
\ncom{\eeq}{\end{equation}}
\ncom{\bea}{\begin{eqnarray*}}
\ncom{\eea}{\end{eqnarray*}}
\ncom{\beqa}{\begin{eqnarray}}
\ncom{\eeqa}{\end{eqnarray}}
\ncom{\nno}{\nonumber}
\ncom{\non}{\nonumber}
\ncom{\ds}{\displaystyle}
\ncom{\half}{\frac{1}{2}}
\ncom{\mbx}{\makebox{.25cm}}
\ncom{\hs}{\mbox{\hspace{.25cm}}}
\ncom{\rar}{\rightarrow}
\ncom{\Rar}{\Rightarrow}
\ncom{\noin}{\noindent}
\ncom{\sz}{\scriptsize}
\ncom{\Sgm}{\Sigma}
\ncom{\psgm}{\sigma^{\prime}}
\ncom{\dt}{\delta}
\ncom{\Dt}{\Delta}
\ncom{\lmd}{\lambda}
\ncom{\Lmd}{\Lambda}
\ncom{\Th}{\Theta}
\ncom{\eps}{\epsilon}
\ncom{\pcc}{\stackrel{P}{>}}
\ncom{\lp}{\stackrel{L_{p}}{>}}
\ncom{\sspan}{{\rm\,span}}
\ncom{\re}{{\rm Re\,}}
\ncom{\im}{{\rm Im\,}}
\ncom{\sgn}{{\rm sgn\,}}
\ncom{\ba}{\begin{array}}
\ncom{\ea}{\end{array}}
\ncom{\integ}[4]{\int_{#1}^{#2}\,{#3}\,d{#4}}
\ncom{\vspan}[1]{{{\rm\,span}\{ #1 \}}}
\ncom{\dm}[1]{ {\displaystyle{#1} } }
\ncom{\ri}[1]{{#1} \index{#1}}

\newtheorem{remark}{\bf Remark}[section]

\newtheorem{example}{Example}[section]

\newtheoremstyle
    {remarkstyle}
    {}
    {11pt}
    {}
    {}
    {\bfseries}
    {:}
    {     }
    {\thmname{#1} \thmnumber{#2} }

\theoremstyle{remarkstyle}

\def \R{{{\rm I{\!}\rm R}}}

\begin{document}
\title
{ Optimal Error Estimates for Semidiscrete Galerkin approximations to the Equations of Motion 
Described by Kelvin-Voigt Viscoelastic Fluid Flow Model}

\author
{Ambit K. Pany{\footnote {Department of Mathematics, Gandhi Institute for Technological Advancement,
Bhubaneswar-752054,India.Email:ambit.pany@gmail.com}},
Saumya Bajpai {\footnote{
Institute of Infrastructure Technology Research and Management, Ahmedabad-380026, India. 
Email: mymsciitm@gmail.com}} and
Amiya K. Pani {\footnote{
 Department of Mathematics, 
Indian Institute of Technology Bombay, Powai, Mumbai-400076, India. Email:akp@math.iitb.ac.in}
}}
\maketitle 
\begin{abstract}
In this paper, the finite element Galerkin method is applied to the equations of motion arising in the 
Kelvin-Voigt viscoelastic fluid flow model, when the forcing function is in $L^{\infty}({\bf L}^2)$.
Some {\it a priori} estimates for the exact solution, which are valid uniformly in time as $t\mapsto \infty$ and 
even uniformly in the retardation time  $\kappa$ an $\kappa\mapsto 0,$ are derived. 
It is  shown that the 
semidiscrete method admits a global attractor.
Further, with the help of {\it a priori} bounds and Sobolev-Stokes projection, optimal error 
estimates for the velocity in $L^{\infty}(\bL^2)$ and $L^{\infty}(\bH^1)$-norms and for the 
 pressure in $L^\infty(L^2)$-norm 
are established.
Since the constants involved in error estimates have an exponential growth in time, therefore, 
in the last part of the article, under certain 
uniqueness condition, the error bounds are established which are valid uniformly in time.  Finally, some
numerical experiments are conducted which confirm our theoretical findings.
\end{abstract}

\noindent{\bf AMS 1991 Classification:}

\section{\normalsize\bf Introduction}
Consider the following system of partial differential equations arising in the Kelvin-Voigt's model
\begin{eqnarray}
\label{1.1}
\frac {\partial \bu}{\partial t}+ \bu\cdot \nabla \bu - \kappa \Delta \bu_t 
-\nu \Delta \bu  + \nabla p
=\f(x,t),\,\,\, x\in \Omega ,\,\,t>0,
\end{eqnarray}
and incompressibility condition
\begin{eqnarray}
\label{1.2}
\nabla \cdot \bu=0,\,\,\,x\in \Omega,\,t>0,
\end{eqnarray}
with initial and boundary conditions
\begin{eqnarray}
\label{1.3}
\bu(x,0)= \bu_0 \;\;\;\mbox {in}\;\Omega,\;\;\;\;\; \bu=0,\;\; \;
\mbox {on}\; \partial \Omega,\; t\ge 0,
\end{eqnarray}
where, $\Omega$ is a bounded convex polygonal or polyhedral  domain in $\R^d, d=2,3$  with boundary $\partial \Omega.$  
Here, $\nu$ is the coefficient of kinematic viscosity and $\kappa$ is the retardation time or the time 
of relaxation of deformations. In the context of viscoelastic fluid, this model was first introduced by Pavlovskii
\cite{P-71}, who called it as a model describing the motion of weakly concentrated water-polymer solutions. 
It was  called Kelvin-Voigt model by  Oskolkov \cite{O-73} and his collaborators. 
Subsequently, Cao {\it et. al.} \cite{CLT-06}  proposed it as  a smooth, inviscid regularization of the 
 2D and 3D-Navier-Stokes equations.
For applications of such models in organic polymer and  food industry, and in the mechanisms of diffuse axonal 
injury, etc., 
we refer to \cite{bs05}, \cite{bs06} and \cite {css02}. 

Earlier, based on the analysis of Ladyzenskaya \cite{L69} in the context of 
 Navier Stokes equations, Oskolkov \cite{APO1}-\cite{APO2} 
have proved  existence of a unique `almost' classical
solution  in finite time interval for the problem
(\ref{1.1})-(\ref{1.3}). Subsequently, further investigations on solvability were continued  
by group members of Oskolkov, see \cite{OR1} and \cite{OR2}. 

On numerical analysis of such problems, Oskolkov {\it et a.} \cite{APO3} have discussed the 
convergence analysis of the spectral Galerkin approximation for all $t \geq 0$ assuming that the exact 
solution is asymptotically stable as $t \rightarrow \infty.$  Subsequently, Pani {\it et a.} \cite{PPDY} have
 applied 
a variant of nonlinear semidiscrete spectral Galerkin method and optimal error estimates are proved. 
It is, further, shown that {\it a priori} error estimates are valid uniformly in time under uniqueness 
assumption.
Recently, Bajpai {\it et al.} \cite{ANS} have applied finite element Galerkin 
methods for the problem (\ref{1.1})-(\ref{1.3}) with the forcing function $\f =0.$ They have proved 
{\it a priori} bounds for the exact solution in $3D$ and established exponential decay property. 
With an introduction of the Sobolev-Stokes projection, they have derived optimal error estimates,
 which again preserve the exponential decay property. In \cite{ANS1}, 
completely discrete schemes which are based on both backward Euler and second order backward 
difference methods are analyzed and optimal error bounds which again preserve exponential decay property
are established. For related articles in the context of Oldroyd 
 viscoelastic model, we refer to \cite{GP11}-\cite {HLSST}, \cite{PY,PYD}, \cite{WHS}-\cite{WLH12}.

In this paper, we,  further, continue the investigation on finite element approximation to the problem 
(\ref{1.1})-(\ref{1.3}) when the non-zero forcing function $\f$ belongs to $L^{\infty}(\bL^2).$  
This is crucial particularly in the study of  the dynamical system (\ref{1.1})-(\ref{1.3}), 
when the forcing function is assumed  to be time independent.
The major results obtained in this paper are summarized as  follows:
{\begin{itemize}
\item [(i)] New regularity results for the solution of (\ref{1.1})-(\ref{1.3}) even in $3D$, which
are valid uniformly in time are derived and as a consequence, existence  of a global attractor is proved.
It is further shown that these estimates hold uniformly in $\kappa$ as $\kappa \mapsto 0.$
\item [(ii)] When $\f$ is independent of time, it is, further, established that the semi-discrete finite element 
method admits a discrete  global  attractor.
\item [(iii)] Based on the Sobolev-Stokes projection introduced earlier in \cite{ANS}, optimal error estimates for the semidiscrete Galerkin approximations to the velocity in $L^{\infty}({\bf L}^2) $-norm as well as 
in $L^{\infty}({\bf H}^1_0)$-norm and to the pressure in $L^{\infty}( L^2)$-norm  are derived
with error bounds depending on exponential in time.
\item [(iv)] Moreover, it is proved under uniqueness assumption that error estimates are valid uniformly in time.
\end{itemize}}
Note that for (i), exponential weight functions in time are used which help us to derive regularity result
for all $t>0.$ A special care is taken to show that these estimates are valid uniformly in $\kappa$ as
$\kappa\mapsto 0.$ When $\f$ is independent of time, based on uniform estimates in time existence of a global 
attractor is shown for the semidiscrete scheme. For (iii), a use of Sobolev-Stokes projection as an 
intermediate projection  helps us to retrieve optimal error estimates for the velocity vector in $L^{\infty}(\bL)$-norm.
When either $\f=0$ or $\f=O(e^{-\alpha_0 t}),$ we derive, as in \cite{ANS}, exponential decay property not only 
for the solution, but also for error estimates.

This paper is organized as follows. In Section {\bf 2},
we discuss the weak formulation and state some basic assumptions. Section {\bf 3} is devoted 
to development of {\it a priori} bounds for the exact solutions. In Section ${\bf 4}$, we describe the 
semidiscrete Galerkin approximations and derive  {\it a priori} estimates with discrete global attractor for
the semidiscrete solutions. In Section {\bf 5}, 
we establish optimal error estimates for the velocity. Section {\bf 6} deals with the optimal error 
estimates for the pressure. In Section {\bf 7}, results of numerical experiments, which confirm our
theoretical estimates, are established. 
\section{ \normalsize \bf Preliminaries and Weak formulation}
\setcounter{tdf}{0}
\setcounter{ldf}{0}
\setcounter{cdf}{0}
\setcounter{equation}{0}
In this section, we define $\mathbb {R}^d,\;( d=2,3)$-valued function
spaces using boldface letters as 
\be
{\bf H}_0^1 = (H_0^1(\Omega))^d, \;\;\; {\bf L}^2 = (L^2(\Omega))^d \;\;
\mbox {and }\;\; {\bf H}^m=(H^m(\Omega))^d,
\ee
where $L^2(\Omega)$ is the space of square integrable functions defined 
in $\Omega $ with inner product $(\phi,\psi) =\displaystyle{\int_{\Omega}}\phi(x) 
\psi(x)\,dx $ and norm  $\|\phi\| = \left(\displaystyle{\int_{\Omega}} 
|\phi(x)|^2\,dx\right)^{1/2}$. Further, $H^m(\Omega)$ denotes the standard 
Hilbert Sobolev space of order $m\in \mathbb{N^+}$ with norm $\|\phi\|_m=
 \left(\displaystyle{\sum_{|\alpha| \leq m}}\displaystyle{\int_{\Omega}} 
|D^{\alpha}\phi|^2\,dx\right)^{1/2}$. Note that $\bH^1_0$ is equipped 
with a norm 
$$\|\nabla\bv\|= \left(\displaystyle{\sum_{i,j=1}^{d}}
(\partial_j v_i, \partial_j v_i)\right)^{1/2}
=\left(\displaystyle{\sum_{i=1}^{d}}
(\nabla v_i, \nabla v_i)\right)^{1/2}.$$
Further, introduce divergence free spaces :
\be
{\bf J}_1 = \{{\bphi} \in {\bf {H}}_0^1 : \nabla \cdot \bphi = 0\}
\ee
and
\be
{\bf J}= \{\bphi \in {\bf {L}}^2 :\nabla \cdot \bphi = 0\;\;
{\mbox {\rm in}}
\;\; \Omega, \;\;\bphi \cdot {\bf {n}} |_{\pr \Omega} = 0\;\;
{\mbox {\rm holds} \;  {\rm weakly}} \},
\ee
where ${\bf {n}}$ is the outward normal to the boundary
$\pr \Omega$ and $\bphi \cdot {\bf {n}} |_{\pr \Omega} = 0$ should
be understood in the sense of trace in $\bH^{-1/2}(\partial \Omega)$,
see \cite{temam}. Let $H^m/{\rm I\!R} $ be the quotient space with norm $\| p\|_{H^m /{\rm I\!R}}
 = \inf_{c\in{\rm I\!R} }\| p+c\|_m$. For a Banach Space $X$ with norm $\|\cdot\|_X,$
 let $L^p(0, T; X)$ denote the space of measurable
 $X$- valued functions $\bphi$ on $ (0,T)$ such that $
\int_0^T \| \bphi (t)\|^p_X dt <  \infty \;\;\;{\mbox {\rm if}}\;\;
1 \le p < \infty$ and for $p=\infty$, ${\dps{ess \sup_{0<t<T}
}} \| \bphi (t)\|_X  < \infty.$
Now, set $P:{\bf L}^2\longrightarrow {\bf J} $ as the ${\bf L}^2$-
orthogonal projection.

Throughout this paper, the following assumptions are made. \\
\noindent
({\bf A1}). 
Setting
$-\tl{\Delta} = -P\Delta:{\bf J}_1 \cap {\bf H}^2
\subset {\bf J} \rightarrow  {\bf J}
$
as the Stokes operator, assume that the following regularity result holds:
\ben \label{2.1a} 
\| \bv\|_2 \le C \| \tl{\Delta} \bv \| \;\;\;\;\forall \bv\in {\bf J}_1
\cap {\bf H}^2.
\een
The above assumption is valid as the domain $\Omega$ is a convex polygon or convex polyhedron.
Note that the following Poincar\'e inequality \cite{HR82} holds true:
\ben \label{2.1*} 
\|\bv\|^2 \le \lambda_1^{-1} \| \nabla \bv \|^2
\;\; \forall \bv \in {\bf H}_0^1(\Omega),
\een
where $ \lambda_1^{-1}$, is the best possible positive constant depending on the domain $\Omega.$
Further, observe that 
\ben \label{2.1**} 
\| \nabla \bv\|^2 \le \lambda_1^{-1} \| \tl{\Delta} \bv\|^2\;\;
\forall \bv  \in {\bf J}_1 \cap {\bf H}^2.
\een

\noindent
({\bf A2}). There exists a positive constant $ M $ such that the initial velocity $\bu_0$ and the external 
force $\f,\f_t$ satisfy for 
$t \in (0,\infty)$ 
\be
\bu_0 \in {\bf H}^2\cap {\bf  J}_1,\, \f,\;\f_t \in L^{\infty}(0,\infty;\,\bL^2)\,\, \rm{with}
\,\, \|\bu_0\|_2 \le M,\,\, 
\displaystyle{\mathop{ess\,\sup}_{0<t< \infty}}\,\|\f(\cdot,t)\|\leq M.\nonumber
\ee
Now, the weak formulation of (\ref{1.1})-(\ref{1.3}) is to seek  
 a pair of functions $(\bu(t), p(t))\in \bH_0^1\times L^2/{\rm I\!R}$  with $\bu(0)= \bu_0$, such that 
 for all $t>0$
\begin{eqnarray}\label{2.2}
  \left.
 \begin{array}{rcl}
 &&(\bu_t, \bphi) +\kappa (\nabla \bu_t, \nabla \bphi) +\nu (\nabla \bu, \nabla \bphi)+ ( \bu \cdot
\nabla \bu, \bphi) = ( p, \nabla \cdot \bphi)+(\bf f, \bphi) \;\forall \bphi \in {\bf H}_0^1,\\
 &&(\nabla \cdot \bu, \chi) = 0 \;\;\; \;\;\;\;\forall \chi \in L^2.
 \end{array}
 \right\}
 \end{eqnarray}
Equivalently,  find  $\bu(t) \in {\bf J}_1 $  with $\bu(0)= \bu_0$  such that for $t>0$
\begin{align}
\label{2.3}
(\bu_t, \bphi) +\kappa (\nabla \bu_t, \nabla \bphi )+\nu (\nabla \bu, \nabla \bphi )
&+( \bu \cdot \nabla \bu, \bphi)
=(\bf f, \bphi)  \;\;\;\forall \bphi \in {\bf J}_1.
 \end{align}
Define the trilinear form $b(\cdot, \cdot, \cdot)$ as 
$$
b(\bv, \bw,\bphi):= \frac{1}{2} (\bv \cdot \nabla \bw , \bphi)
- \frac{1}{2} (\bv \cdot \nabla \bphi, \bw),\;\;\bv, \bw, \bphi \in \bH_0^1.
$$
Note  for $\bv \in \bJ_1,$  $\bw, \bphi \in \bH^1_0$ that
$ 
b(\bv, \bw,\bphi)= (\bv \cdot \nabla \bw , \bphi).
$
Because of antisymmetric property of the trilinear form, it is easy to check that for $,$
\begin{equation} \label{ASP}
b(\bv,\bw,\bw) =0\;\;\;\forall \bv,\bw\in \bJ_1.
\end{equation}

\section{ \normalsize \bf A priori estimates for the exact solution}
\setcounter{tdf}{0}
\setcounter{ldf}{0}
\setcounter{cdf}{0}
\setcounter{equation}{0}
In this section, some {\it a priori} bounds for the solution ${(\bu, p)}$ of (\ref{2.2}) are derived. 
Since these results differ from \cite{ANS} in the sense that $0\neq\f\in L^{\infty}(\bL^2)$ in the present article, therefore, 
only the major differences in the analysis are indicated. 
\noindent
\begin{ldf}\label{L42}
Let the assumptions {\rm{({\bf A1})}}-{\rm{({\bf A2})}} hold true, and let 
$\displaystyle{0 <\alpha<\frac{\nu\lambda_1}{4\left(1+\kappa\lambda_1 \right)}}$. Then, 
the solution $\bu$ of (\ref{2.3}) satisfies for all $t > 0$
\begin{align}
 \label{nc6}
\Big(\|\bu(t)\|^2 &+\kappa\|\nabla \bu(t)\|^2\Big) + \beta e^{-2\alpha t}
\int_0^t{e^{2\alpha s}\|\nabla \bu(s)\|^2}\,ds\nonumber \\
&\le  e^{-2\alpha t}(\|\bu_0\|^2+\kappa \|\nabla\bu_0\|^2)+ \left(\frac{1-e^{-2\alpha t}}{2\nu\lambda_1\alpha}\right)
\|{\bf f}\|^2_{L^{\infty}(\bL^2)}=:K_0(t) \nonumber\\
&\le (\|\bu_0\|^2+\kappa \|\nabla\bu_0\|^2)+ \left(\frac{1}{2\nu\lambda_1\alpha}\right)
\|{\bf f}\|^2_{L^{\infty}(\bL^2)}=:  K_{0,\infty},\;\;\;\;t>0.
\end{align}
where $\beta= \nu-2\alpha (\kappa+\lambda_1^{-1})\geq \nu/2>0,$ and $K_{0, \infty}= \displaystyle{\mathop{\sup}_{t \in [0,\infty)}}\, K_0(t)$. Moreover,
 \begin{eqnarray}
  \label{nc1}
  \mathop{\lim \,\sup}_{t\rightarrow \infty}  \|\nabla\bu(t)\|\leq \left( \frac{1}{\lambda_1\nu^2}\right)\|\f\|_{L^{\infty}{(0,\infty;\;\bL^2)}}. 
 \end{eqnarray}
\end{ldf}
\noindent{\it Proof.} 
Set $\hat{\bu}(t)=e^{\alpha t}\bu(t)$ for some $\alpha > 0$ in (\ref{2.3}). 
Then, choose $\bphi=\hat\bu$ in (\ref{E42}) and use (\ref{ASP}) in the resulting equation to arrive at
\begin{eqnarray}
\label{ne1}
\frac{1}{2}\frac{d}{dt}(\|\hat \bu\|^2+\kappa \|\nabla\hat\bu\|^2)+\left(\nu-\alpha (\kappa+\lambda_1^{-1})\right)
\|\nabla\hat \bu\|^2 \leq (\hat{\bf f},\hat \bu).
\end{eqnarray}
Now, estimate the right-hand side of (\ref{ne1}) as
\begin{eqnarray}
\label{2.6}
|(\hat {\bf f},\hat\bu)|\le \|\hat {\bf f}\| \|\hat\bu\| \le \frac{1}{\sqrt\lambda_1}\|\hat {\bf f}\| \|\nabla\hat\bu\| 
\le \frac{\nu}{2}\|\nabla\hat\bu\|^2 + \frac{1}{2\nu\lambda_1}\|\hat {\bf f}\|^2.
\end{eqnarray}
Substitute (\ref{2.6}) in (\ref{ne1}), use kickback argument and $\beta= \nu-2\alpha (\kappa+\lambda_1^{-1})
=\nu/2- (\nu/2-2\alpha (\kappa+\lambda_1^{-1})) \geq \nu/2>0$ to obtain
\begin{eqnarray}
 \label{2.8}
\frac{d}{dt}(\|\hat \bu\|^2+\kappa \|\nabla\hat\bu\|^2)+\beta \|\nabla\hat \bu\|^2 \le \frac{1}{\nu\lambda_1}
\|\hat {\bf f}\|^2.
\end{eqnarray}
Integrate with respect to time from 0 to $t$, then multiply by $ e^{-2\alpha t}$ and use the assumption {\rm{({\bf {A2}})}} 
as well as the fact that 
\begin{eqnarray}
\label{2.81}
e^{-2 \alpha t} \int_0^t e^{2 \alpha s} ds = \frac{1}{2\alpha} (1 - e^{-2 \alpha t})
\end{eqnarray}
to complete the proof of (\ref{nc6}). 

Note that the second term  on the left had side of  (\ref{nc6}) is nonnegative 
and hence, it can be dropped. Then taking limit superior  as $t\longrightarrow \infty$ for the remaining terms
on both sides, we arrive at   
\begin{equation}\label{sup-time}
\mathop{\lim \,\sup}_{t\longrightarrow \infty}
(\|\bu(t)\|^2
+\kappa\|\nabla\bu(t)\|^2) \leq \left(\frac{1}{2\nu\lambda_1\alpha}\right)
\|{\bf f}\|^2_{L^{\infty}(\bL^2)}.
\end{equation}
For (\ref{nc1}), we
rewrite (\ref{ne1}) as :
\begin{eqnarray}
\frac{1}{2}\frac{d}{dt}(\|\hat \bu\|^2+\kappa \|\nabla\hat\bu\|^2)+\nu\|\nabla\hat \bu\|^2 \leq (\hat{\bf f},\hat \bu)
+\alpha (\|\hat\bu\|^2+\kappa\|\nabla\hat\bu\|^2).\nonumber
\end{eqnarray}
Integrate with respect to time and then, divide the resulting equation by $e^{-2\alpha t}$ to arrive at
\begin{align}\label{nc4}
\left(\|\bu(t)\|^2+\kappa \|\nabla\bu(t)\|^2\right)&+\nu e^{-2\alpha t}\int_0^t e^{2\alpha s}\|\nabla\bu(s)\|^2ds\leq e^{-2\alpha t}
(\|\bu_0\|^2+\kappa\|\nabla\bu_0\|^2)\nonumber\\
+\frac{\|\f\|^2_{L^{\infty}(\bL^2)}}{2\alpha\lambda_1\nu}& (1-e^{-2\alpha t})+2 \alpha e^{-2\alpha t}\int_0^te^{2\alpha s}(\|\bu(s)\|^2
+\kappa\|\nabla\bu(s)\|^2)ds.
\end{align}
Now, the first term on the left hand side of (\ref{nc4}) is nonnegative which can then be dropped.
Taking limit superior on the both sides of  (\ref{nc4}) for the remaining terms and using L' Hospital rule, we note
that
\begin{align}
\label{nc5*}
 \mathop{\lim \,\sup}_{t\longrightarrow \infty}2 \alpha e^{-2\alpha t}\int_0^te^{2\alpha s}(\|\bu(s)\|^2
+\kappa\|\nabla\bu(s)\|^2)ds
&=\mathop{\lim \,\sup}_{t\longrightarrow \infty}
(\|\bu(t)\|^2
+\kappa\|\nabla\bu(t)\|^2),
\end{align}
\begin{align}
\label{nc6*}
 \mathop{\lim \,\sup}_{t\longrightarrow \infty}\nu e^{-2\alpha t}\int_0^t e^{2\alpha s}\|\nabla\bu(s)\|^2ds
 &=\frac{\nu}{2\alpha}\mathop{\lim \,\sup}_{t\longrightarrow \infty}\|\nabla\bu(t)\|^2,
\end{align}
and hence, using (\ref{sup-time}) we arrive at
\begin{align}
 \mathop{\lim \,\sup}_{t\longrightarrow \infty}\|\nabla\bu(t)\|\leq \left( \frac{1}{\lambda_1\nu^2}\right) 
 \|{\bf f}\|_{L^{\infty}(0,\infty;\;\bL^2)}.\nonumber
\end{align}
This completes the rest of the proof. \hfill {$\Box$}

\begin{remark}
As a consequence of Lemma \ref{L42}, we obtain from (\ref{2.8}) with $\alpha=0$ the following 
estimate
\begin{eqnarray}
 \label{2.8a}
\frac{d}{dt}(\|\bu\|^2+\kappa \|\nabla \bu\|^2)+\nu \|\nabla  \bu\|^2 \le \frac{1}{\nu\lambda_1}
\| {\bf f}\|^2.
\end{eqnarray}
On integration with respect to time from $t$ to $t+T_0$, and using (\ref{nc6}) of Lemma \ref{L42}, we 
obtain for fixed $T_0>0$ and $t\geq 0$
\begin{eqnarray}
 \label{nc6a}
\nu \int_{t}^{t+T_0} \|\nabla  \bu\|^2 \,ds &\le& K_0(t)+\frac{T_0}{\nu\lambda_1}
\| {\bf f}\|^2 \nonumber\\
&\le& K_{0,\infty} + \frac{T_0}{\nu\lambda_1}
\| {\bf f}\|^2.
\end{eqnarray}
Taking limit superior on both sides of (\ref{nc6a}), we now arrive at
\begin{eqnarray}
 \label{nc6b}
\nu  \mathop{\lim \,\sup}_{t\longrightarrow \infty}\;\int_{t}^{t+T_0} \|\nabla  \bu\|^2 \,ds 
\le K_{0,\infty} + \frac{T_0}{\nu\lambda_1}
\| {\bf f}\|^2.
\end{eqnarray}
\end{remark}
\begin{remark}
Note that if $\f \in L^{\infty}(\bH^{-1}),$ where ${\bH^{-1}}$ is the topological dual of 
$ {\bH^1_0},$ then following the proof of the Lemma \ref{L42}, obtain 
\begin{align}
 \label{nc6-1}
\|\bu(t)\|^2 &+\kappa\|\nabla \bu(t)\|^2+ \beta e^{-2\alpha t}
\int_0^t{e^{2\alpha s}\|\nabla \bu(s)\|^2}\,ds\nonumber \\
&\le  e^{-2\alpha t}(\|\bu_0\|^2+\kappa \|\nabla\bu_0\|^2)+ \left(\frac{1-e^{-2\alpha t}}{2\nu\alpha}\right)
\|{\bf f}\|^2_{L^{\infty}({(\bH^1_0)^{*}})}=K^*_0(t) \le K^*_{0,\infty} ,\;\;\;\;t>0.
\end{align}
\end{remark}
\begin{remark}
Earlier, Oskolkov \cite{APO2} has proved the existence of a unique weak solution to the problem (\ref{1.1})-
(\ref{1.3}) for finite time, but 
the proof can not be extended  to all $t>0$ as the constants involved in {\it a priori} estimates depend on 
exponentially in time.  Now, 
 using Bubnov Galerkin method with {\it a priori} bounds in Lemma \ref{L42} and standard weak 
 compactness arguments, it can be shown that there exists a unique global weak solution $\bu$ to the problem 
(\ref{2.3}) for all $t>0$. Further, it is easy to check that the 
 problem (\ref{2.3}) generates a continuous semigroup $S(t):\bJ_1\rightarrow\bJ_1,\,\,t\in [0,\infty)$. 
Therefore, the result of \cite{KT} shows that if $\f \in L^{\infty}({\bH^{-1}}),$ then the semigroup
$S(t)$ has an absorbing ball  
$$B_{\brho}(0):\{\bv \in \bJ_1 : \Big(\|\bv\|^2 +\kappa \|\nabla \bv\|^2\Big)^{1/2} \leq \rho \}$$ 
with $\rho$ given by 
 $$\brho^2=\left(\frac{1}{\alpha \nu }\right)\|\f\|^2_{L^{\infty}((\bH^1_0)^{*})}.
 $$
 Hence, it may be easily shown that the problem has  a global attractor ${\mathcal{A}}_1 \subset \bJ_1.$ 
\end{remark}

\begin{ldf}\label{L43}
Let 
assumptions {\rm{({\bf A1})-({\bf A2})}} hold true. Then, for $\displaystyle{0 <\alpha< \frac{\nu \lambda_1}{4\left(1+\lambda_1\kappa\right)}}$  and for all $t>0$
\begin{eqnarray*}
\|\nabla\bu(t)\|^2 &+&\kappa\|\tilde \Delta \bu(t)\|^2 + \beta e^{-2\alpha t}\int_{0}^t  e^{2\alpha s}
\|\tilde \Delta \bu(s)\|^2 \,d s 
 \le e^{-2\alpha t}(\|\nabla\bu(0)\|^2 +\kappa\|\tilde \Delta \bu(0)\|^2) \nonumber\\
&+& C(\nu,\alpha)\left(\frac{K_{0, \infty}^{\ell+2}}{\kappa^{\ell}}(1-e^{-2\alpha t})+(1-e^{-2\alpha t})\|
{\bf f}\|^2_{L^{\infty}(\bL^2)}\right)
=K_1(t)\leq K_{1,\infty}
\end{eqnarray*}
holds, where $\beta= \nu-2\alpha (\kappa+\lambda_1^{-1})\geq \nu/2 >0,$  for $d=2$, $\ell=1,$ and when
$d=3,$ $\ell=3.$
\end{ldf}
\noindent{\it Proof.} 
Set $\hat\bu=e^{\alpha t}\bu$ and use the definition of the Stokes operator $\tilde\Delta$ to 
rewrite (\ref{E42}) as
\begin{align}\label{E43}
(\hat\bu)_t &-\alpha\; \hat\bu -\kappa\;\tilde\Delta \hat\bu_t +\kappa\alpha\;\tilde\Delta 
\hat\bu -\nu\;\tilde\Delta\hat\bu =-e^{-\alpha t} (\hat\bu\cdot\nabla \hat{\bu})+ \hat{\bf f}\;\;
\,\,\,\forall \bphi\in{\bf J}_1.
\end{align}
\noindent
Multiply (\ref{E43}) by $-\tilde\Delta\hat\bu$  and integrate over $\Omega.$ A use of integration by parts
with  (\ref{2.1*}) and $-(\hat\bu_t,\tilde\Delta\hat\bu)
=\displaystyle{\frac{1}{2}\frac{d}{ dt}\|\nabla \hat\bu\|^2}$ 
 leads to 
\begin{align}\label{E44}
\frac{1}{2}\frac{d}{ dt}(\|\nabla \hat\bu\|^2+\kappa\|\tilde\Delta \hat\bu\|^2)+\left(\nu-\alpha (\kappa+\lambda_1^{-1})\right)
\|\tilde\Delta\hat\bu\|^2&=e^{-\alpha t}(\hat\bu\cdot\nabla
\hat\bu,\tilde\Delta \hat\bu)+(\hat{\bf f},-\tilde\Delta\hat\bu)\nonumber\\
&=I_1+I_2.
\end{align}
\noindent
For $I_1,$ we note by generalized H\''{o}lder's inequality that
\begin{equation}
\label {E43-0}
|I_1|\leq e^{-\alpha t } \|\hat \bu\|_{L^4}\;\|\nabla \hat\bu \|_{L^4}\;\|\tilde \Delta\hat \bu\|.
\end{equation} 
When $d=2,$ a use of  Ladyzhenskaya's inequality:
$$\|\hat \bu\|_{L^4} \leq C\; \|\hat \bu\|^{\frac{1}{2}}\;\|\nabla \hat\bu \|^{\frac{1}{2}}\;\;\;
\mbox{and}\;\; \|\nabla \hat\bu \|_{L^4} \leq \|\nabla \hat\bu \|^{\frac{1}{2}}\;
\|\Delta \hat\bu \|^{\frac{1}{2}}.
$$
in (\ref{E43-0}) with the Young's inequality with $p=4$, $q=\frac{4}{3}$, $\epsilon=\frac{2\nu}{9}$ yields
\begin{align}
\label {E43a}
|I_1|& \leq Ce^{-\alpha t } \|\hat \bu\|^{\frac{1}{2}}\|\nabla \hat\bu \|\|\tilde \Delta\hat \bu\|^{\frac{3}{2}}
     \leq C  {\bigg(\frac{1}{\nu}\bigg)}^3 e^{2\alpha t }  \| \bu\|^2 \| \nabla \bu\|^{4} 
     +\frac{\nu}{6} \| \tilde  \Delta \hat \bu \|^{2}. 
\end{align}
When $d=3,$ a use of  Ladyzhenskaya's inequality:
\begin{equation} \label{LI-d3}
\|\hat \bu\|_{L^4} \leq C\; \|\hat \bu\|^{\frac{1}{4}}\;\|\nabla \hat\bu \|^{\frac{3}{4}}\;\;\;
\mbox{and}\;\; \|\nabla \hat\bu \|_{L^4} \leq \|\nabla \hat\bu \|^{\frac{1}{4}}\;
\|\Delta \hat\bu \|^{\frac{3}{4}}.
\end{equation}
in (\ref{E43-0}) with the Young's inequality with $p=8/7$, $q=8$, $\epsilon^p=\frac{4\nu}{21}$ shows
\begin{align}
\label {E43a-1}
|I_1|& \leq Ce^{-\alpha t } \|\hat \bu\|^{\frac{1}{4}}\|\nabla \hat\bu \|\;\|\tilde \Delta\hat \bu\|^{\frac{7}{4}}
     \leq C  {\bigg(\frac{1}{\nu}\bigg)}^7 e^{2\alpha t }  \| \bu\|^2 \| \nabla \bu\|^{8} 
     +\frac{\nu}{6} \| \tilde  \Delta \hat \bu \|^{2}. 
\end{align}
For $I_2,$ an application of  the Cauchy-Schwarz inequality with the Young's inequality leads to 
\begin{eqnarray}
\label{2.9}
|I_2|=|(\hat {\bf f},-\tilde\Delta\hat\bu)|\le \|\hat {\bf f}\| \|\tilde\Delta\hat\bu\| 
\le \frac{\nu}{3}\|\tilde\Delta\hat\bu\|^2 + \frac{3}{2\nu}\|\hat {\bf f}\|^2.
\end{eqnarray}  
 Substitute (\ref{E43a}) and (\ref{2.9})  in (\ref{E44}) to find at
\begin{align}\label{E43e}
 \frac{d}{dt}\left(\|\nabla\hat\bu \|^2+\kappa \|\tilde \Delta \hat\bu\|^2\right)+\left(\nu - 2\alpha(\kappa+\lambda^{-1})\right)
\|\tilde \Delta \hat\bu\|^2 
 \leq C(\nu)\Big( e^{2\alpha t}\| \bu\|^2\;\|\nabla  \bu\|^{2(\ell+1)}+ \|\hat {\bf f}\|^2\Big), 
\end{align}
where $\ell=1,$ when $d=2$ and for $d=3$, $\ell=3.$
Integrate (\ref{E43e}) with respect to time from $0$ to $t$. Then, use Lemma \ref{L42} and $\beta= 
\nu - 2\alpha(\kappa+\lambda^{-1})\geq \nu/2>0 $ 
to arrive at 
\begin{eqnarray} \label{eq:grad}
\|\nabla\bu(t)\|^2&+&\kappa\|\tilde \Delta \bu(t)\|^2 +\beta e^{-2\alpha t}\int_{0}^t  e^{2\alpha s}
\|\tilde \Delta \bu(s)\|^2 \,ds \le e^{-2\alpha t}(\|\nabla\bu_0\|^2 +\kappa\|\tilde \Delta \bu(0)\|^2)\nonumber\\
 &+&C(\nu)e^{-2\alpha t}\int_0^t e^{2\alpha s}\|\bu(s)\|^2\;\|\nabla \bu(s)\|^2\;\|\nabla \bu(s)\|^{2\ell} \; ds
 \nonumber\\
&+& C(\nu,\alpha)(1-e^{-2\alpha t})\|\bf f\|^2_{L^{\infty}(\bL^2)}.
\end{eqnarray}
For the second term one the right hand side of (\ref{eq:grad}), apply Lemma \ref{L42} to obtain
\begin{eqnarray}
\|\nabla\bu(t)\|^2 &+&\kappa\|\tilde \Delta \bu(t)\|^2 +\beta e^{-2\alpha t}\int_{0}^t  e^{2\alpha s}
\|\tilde \Delta \bu(s)\|^2 \,ds 
\leq e^{-2\alpha t}(\|\nabla\bu_0\|^2 +\kappa\|\tilde \Delta \bu_0\|^2)\nonumber\\
&+&C(\nu,\alpha)\left(\frac{K_{0,\infty}^{\ell+2}}{\kappa^{\ell}}(1-e^{-2\alpha t})+\|{\bf f}\|_{L^{\infty}{(\bL^{2})}}^2(1-e^{-2\alpha t})\right).\nonumber
\end{eqnarray} 
This completes the rest of the proof.\hfill {$\Box$}

Note that results in Lemma \ref{L43} are valid uniformly in time  for both $2D$ and $3D$ problems. However, 
constants in those bounds   depend on $1/\kappa,$ which blow up as $\kappa$ tends to zero. Therefore, in the  
following Lemma, we propose to discuss  results which are valid for  all time, but their bounds are
independent of $1/\kappa.$

\begin{ldf}\label{L43-1}
Let 
assumptions {\rm{({\bf A1})-({\bf A2})}} hold true. Then, there exists a positive constant 
$K_{12}=K_{12}(\nu,\alpha, \lambda_1, M)$
such that for $\displaystyle{0 <\alpha< \frac{\nu \lambda_1}{4\left(1+\lambda_1\kappa\right)}}$  
and for all $t>0,$
\begin{eqnarray}\label{dirichlet}
\|\nabla\bu(t)\|^2 +\kappa\|\tilde \Delta \bu(t)\|^2 + \beta e^{-2\alpha t}\int_{0}^t  e^{2\alpha s}
\|\tilde \Delta \bu(s)\|^2 \,d s 
\leq K_{12},
\end{eqnarray}
where $\beta= \nu-2\alpha (\kappa+\lambda_1^{-1})\geq \nu/2 >0.$ For $d=3$, the estimate (\ref{dirichlet}) holds
true under smallness assumption on $M,$ that is, on the data.
\end{ldf}
\noindent{\it Proof.} When $d=2,$ we note from (\ref{eq:grad}) that
\begin{align}\label{n20}
\|\nabla\hat\bu(t)\|^2 &+\kappa\|\tilde \Delta \hat\bu(t)\|^2 +\beta \int_{0}^t  e^{2\alpha s}
\|\tilde \Delta \bu(s)\|^2 \,ds \le (\|\nabla\bu_0\|^2 +\kappa\|\tilde \Delta \bu_0\|^2)\nonumber\\
&+C(\nu)\displaystyle{\int_0^t}\|\hat\f(s)\|^2\;ds
+C(\nu)\int_0^t\|\bu(s)\|^2\|\nabla \bu(s)\|^2\|\nabla \hat\bu(s)\|^2 ds.
\end{align}  
An application of Gronwall's lemma leads to 
\begin{align}\label{n21}
\|\nabla\hat\bu(t)\|^2 &+\kappa\|\tilde \Delta \hat\bu(t)\|^2 +\beta \int_{0}^t  e^{2\alpha s}
\|\tilde \Delta \bu(s)\|^2 \,d s \le \{(\|\nabla\bu(0)\|^2 +\kappa\|\tilde \Delta \bu(0)\|^2)\nonumber\\
&+C(\nu)\displaystyle{\int_0^t}\|\hat\f(s)\|^2ds\}\times exp\left(C(\nu)\int_0^t\|\bu(s)\|^2\|\nabla \bu(s)\|^2 ds\right).
\end{align} 
Apply assumption {\bf (A2)} in (\ref{n21}) to obtain
\begin{align}\label{n22}
\|\nabla\bu(t)\|^2 &+\kappa\|\tilde \Delta \bu(t)\|^2 +\beta \int_{0}^t  e^{2\alpha s}
\|\tilde \Delta \bu(s)\|^2 \,d s \le C(\nu,\alpha,K_{0,\infty})\;exp\left(C(\nu)\int_0^t\|\bu(s)\|^2\|\nabla \bu(s)\|^2 ds\right).
\end{align} 
A use of estimate (\ref{nc6}) of Lemma \ref{L42} with estimate (\ref{nc6b}) in (\ref{n22}) shows 
that for all finite but fixed $0<T_0$ with $0<t\leq T_0$ and for $d=2$
\begin{align}\label{n23}
\|\nabla\bu(t)\|^2 &+\kappa\|\tilde \Delta \bu(t)\|^2 +\beta \int_{0}^t  e^{2\alpha s}
\|\tilde \Delta \bu(s)\|^2 \,d s \le C(\nu,\alpha,K_{0,\infty},T_0).
\end{align}
Since the inequality (\ref{n23}) is valid for all finite, but fixed $T_0,$ now a use of
the following result (\ref{nc1}) from Lemma \ref{L42}
$$\mathop{\lim \,\sup}_{t\rightarrow \infty}\|\nabla\bu\|\leq C$$ 
leads to the boundedness of $\|\nabla\bu(t)\|$ for all $t>0$.
This completes the the proof for $d=2.$

When $d=3,$ that is, the problem in $3D,$ we observe from (\ref{eq:grad}) with $\ell=3$ after multiplying with
$e^{-2\alpha t}$ both sides and using (\ref{nc6}) that
\begin{align}\label{n20a}
\|\nabla \bu(t)\|^2 &+\kappa\|\tilde \Delta \bu(t)\|^2 +\beta e^{-2\alpha t} \int_{0}^t  e^{2\alpha s}
\|\tilde \Delta \bu(s)\|^2 \,ds \le e^{-2\alpha t}(\|\nabla\bu_0\|^2 +\kappa\|\tilde \Delta \bu_0\|^2)\nonumber\\
&+\frac{3}{\nu} e^{-2\alpha t}\displaystyle{\int_0^t}\|\hat\f(s)\|^2\;ds
+C(\nu) e^{-2\alpha t}\;\int_0^t\|\bu(s)\|^2 \;\|\nabla \bu(s)\|^8\; ds\nonumber\\
&\leq  C_1(K_{0,\infty}) + C_2(K_{0,\infty}) \;\int_0^t\;\|\nabla \bu(s)\|^8\; ds
\end{align}  
Setting $\Psi= \|\nabla \bu(t)\|^2 $ and dropping the last two terms on the left hand side of (\ref{n20a})
as these are nonnegative, then we arrive at  
\begin{equation}\label{n20b}
\Psi(t) \leq  C_1(K_{0,\infty}) + C_2(K_{0,\infty}) \;\int_0^t\;\Psi^4(s)\; ds
\end{equation}  
This integral inequality holds true for all finite time $t>0$ provided both $C_1(K_{0,\infty})$ and 
$C_2(K_{0,\infty})$ are sufficiently small, that is, under the assumption that the condition $(A_2)$ is valid
for sufficiently small $M.$ Therefore, the boundedness of $\|\nabla \bu(t)\|$ is proved for all finite, but fixed
$t>0$ and for sufficiently smallness assumption on both initial data and forcing function. The rest of the 
analysis follows as in $2D$ case, that is, when $d=2,$ using the estimate (\ref{nc1}). This completes 
the rest of the proof.
\hfill{$\Box$}

\begin{ldf}
\label{L35}
Under assumptions {\rm{({\bf A1})-({\bf A2})}}, there exists a constant 
$C=C(\nu,\alpha, \lambda_1, M)$ such that the following holds true for 
$\displaystyle{0 < \alpha< \frac{\nu \lambda_1}{4\left(1+\lambda_1\kappa\right)}}$ and for all $t>0$
\begin{eqnarray*}
 e^{-2\alpha t}\int_0^t e^{2\alpha s}( \|\bu_{t}(s)\|^2+2\kappa\| \nabla\bu_{t}(s) \|^2)\;ds 
 +\nu\|\nabla\bu(t)\|^2\leq C.\nonumber
\end{eqnarray*}
\end{ldf}
\noindent
 {\it Proof.} 
Choose  $\bphi=e^{2\alpha t}\bu_t$ in (\ref{2.3}) to arrive at  
\begin{eqnarray}
\label{L35.1*}
 e^{2\alpha t}(\|\bu_t\|^2+\kappa\|\nabla\bu_t\|^2)+\frac{\nu}{2} e^{2\alpha t} \frac{d}{dt}\|\nabla \bu\|^2 
=e^{2\alpha t}({\bf f}, \bu_t)-e^{2\alpha t}(\bu.\nabla\bu,\bu_t).
\end{eqnarray}
For the nonlinear term on the right hand side of (\ref{L35.1*}), use  Sobolev imbedding theorem to obtain
\begin{eqnarray}
\label{L35.1**}
|(\bu.\nabla\bu,\bu_t)|\leq C \|\bu\|_{\bL^4}\;\|\nabla \bu\|_{\bL^4}\; \|\bu_t\|
\leq C \|\nabla \bu \|\; \|\tilde{\Delta}\bu\| \; \|\bu_t\|. 
\end{eqnarray}
Use (\ref{L35.1**}) in (\ref{L35.1*}), then  integrate the resulting inequality with respect to time 
from $0$ to $t$ and apply the Young's inequality. Then, multiply 
the resulting equation by $e^{-2\alpha t}$ to arrive at
\begin{align}
e^{-2\alpha t}\int_0^t e^{2\alpha s}(\|\bu_{t}(s)\|^2&+2\kappa\| \nabla\bu_{t}(s) \|^2)ds
+\nu \|\nabla\bu(t)\|^2 \leq Ce^{-2\alpha t}\|\nabla \bu_0\|^2+e^{-2\alpha t}\int_0^t e^{2\alpha s}\|\nabla \bu(s)\|^2 ds\nonumber\\
&+ e^{-2\alpha t}\int_0^t e^{2\alpha s}\|{\bf f}(s)\|^2 ds+
e^{-2\alpha t}\int_0^t e^{2\alpha s}\|\nabla \bu(s)\|^2\|\tilde \Delta \bu(s)\|^2 ds.
\end{align}
A use of Lemmas \ref{L42} with {\ref{L43-1}} leads to the desired result and this concludes the proof.
 \hfill {$\Box$}
\begin{ldf}
\label{L34}
Let the assumptions {\rm{({\bf A1})-({\bf A2})}} hold true. Then, there exists a positive constant 
$C=C(\nu, \alpha, \lambda_1, M)$
such that for all $t>0$
\begin{eqnarray*}
 \|\bu_t(t)\|^2+\kappa \|\nabla \bu_t(t)\|^2+\nu e^{-2\alpha t}\displaystyle{\int_0^t}e^{2\alpha s}\|\nabla \bu_t(s)\|^2ds\leq C.
\end{eqnarray*}
\end{ldf}
\noindent{\it Proof.} Differentiate (\ref{2.3}) with respect to time to obtain
\begin{align}\label{n11}
(\bu_{tt}, \bphi) +\kappa (\nabla \bu_{tt}, \nabla \bphi )+\nu (\nabla \bu_t, \nabla \bphi )
&=-(\bu_t \cdot \nabla \bu, \bphi)-(\bu \cdot \nabla \bu_t, \bphi)
+(\bf f, \bphi)  \;\;\;\forall \bphi \in {\bf J}_1.
\end{align}
Choose $\bphi=\bu_t$ in (\ref{n11}) with $(\bu\cdot \nabla \bu_t, \bu_t)=0$ to find that
\begin{eqnarray}
\label{L34.1}
 \frac{1}{2}\frac{d}{dt}(\|\bu_t\|^2+\kappa\|\nabla \bu_t\|^2)+\nu\|\nabla\bu_t\|^2=
 -(\bu_t \cdot \nabla \bu, \bu_t)
+(\bf f,  \bu_t).
\end{eqnarray}
\noindent
Apply 
the Ladyzenskaya's inequality (\ref{LI-d3}) for $d=3$  and the Young's inequality 
(with $p=8$ and $q=8/7$) to arrive at
\begin{align}\label{n12}
 (\bu_t \cdot \nabla \bu, \bu_t)
 &\leq C \|\bu_t\|^{1/4}\|\nabla\bu\|\|\nabla\bu_t\|^{7/4}\nonumber\\
 &\leq C(\nu)\;\|\nabla\bu\|^8\;\|\bu_t\|^2+\frac{\nu}{4}\;\|\nabla\bu_t\|^2.
\end{align}
A use of the Cauchy-Schwarz inequality with the Young's inequality leads to 
\begin{eqnarray}
\label{L34.2}
  (\f,\bu_t)\leq \|\f\|\;\|\bu_t\|\leq \frac{1}{\sqrt{\lambda_1}} \|\f\|\;\|\nabla \bu_t\| \leq 
  \frac{1}{\lambda_1\;\nu}\|\f\|^2+\frac{\nu}{4}\;\|\nabla\bu_t\|^2.
 \end{eqnarray}
Substitute (\ref{n12})-(\ref{L34.2}) in (\ref{L34.1}) and then  multiply by $e^{2\alpha t}.$ An 
application of {\it a priori} estimates from 
Lemma \ref{L43-1}, \ref{L35} yields
\begin{align}\label{n14}
 \frac{d}{dt}e^{2\alpha t}(\|\bu_t\|^2+\kappa\|\nabla \bu_t\|^2)+\nu e^{2\alpha t}
 \|\nabla\bu_t\|^2 &\leq C(\nu,\lambda_1)e^{2\alpha t}(\|\bu_t\|^2+\|\f\|^2)\nonumber\\
 &+2\alpha e^{2\alpha t}(\|\bu_t\|^2+\kappa\|\nabla \bu_t\|^2).
\end{align}
Integrate (\ref{n14}) from $0$ to $t$ with respect to time to obtain
\begin{align}\label{n15}
 &\|\bu_t\|^2+\kappa\|\nabla \bu_t\|^2+\nu e^{-2\alpha t}\displaystyle{\int_0^t}e^{2\alpha s}
 \|\nabla\bu_t(s)\|^2\;ds\leq e^{-2\alpha t}(\|\bu_t(0)\|^2+\kappa\|\nabla \bu_t(0)\|^2)\nonumber\\
 &+Ce^{-2\alpha t}\displaystyle{\int_0^t}e^{2\alpha s}(\|\bu_t(s)\|^2+\|\f(s)\|^2)ds
 +2\alpha e^{-2\alpha t}\displaystyle{\int_0^t} e^{2\alpha s}(\|\bu_t(s)\|^2+\kappa\|\nabla \bu_t(s)\|^2)ds.
\end{align}
From (\ref{2.3}), it may be observed that
\begin{align}\label{n16}
 \|\bu_t\|^2+\kappa\|\nabla\bu_t\|^2&\leq C(\|\tilde\Delta \bu\|^2+\|\f\|^2+\|\bu\|^2\|\nabla\bu\|^4)\nonumber\\
 &\leq C(\lambda_1)(\|\tilde\Delta \bu\|^2+\|\f\|^2).
\end{align}
Using (\ref{n16}) (see, the proof in \cite{HR82} pp 285, eq (2.19)), we can define (\ref{n16}) at $t=0$. 
A use of Lemma \ref{L35} with {\bf (A2)} and (\ref{n16}) in (\ref{n15}) establishes the desired estimates.
This completes the rest of the proof \hfill {$\Box$} 


\begin{ldf}\label{L37}
 Let assumptions {\rm{({\bf A1})-({\bf A2})}} hold. Then, there exists
a positive constant $C=C(\nu,\alpha,\lambda_1, M)$
such that for $\displaystyle{0<\alpha< \frac{\nu \lambda_1}{4\left(1+\lambda_1\kappa\right)}}$  and 
for all $t>0$,
\begin{eqnarray}\label{estimate-Delta}
\nu \|\tilde\Delta\bu(t)\|^2 + e^{-2\alpha t}\int_0^t e^{2\alpha s}( \|\nabla\bu_{t}(s)\|^2
+\kappa\| \tilde\Delta\bu_{t} (s)\|^2)ds \leq C.
\end{eqnarray}
Moreover, the following estimate hold:
\begin{eqnarray}\label{estimate-Delta-1}
\kappa \|\tilde\Delta\bu_t(t)\| \leq C.
\end{eqnarray}
\end{ldf}
\noindent{\it Proof.} Rewrite (\ref{2.3}) as 
\begin{eqnarray}
\label{L35.1}
 \bu_t -\kappa\;\tilde \Delta\bu_t -\nu\;\tilde\Delta \bu + \bu\cdot\nabla\bu 
 = \bf f \,\,\, \,\,\,\forall \bphi\in\bJ_1.
\end{eqnarray}
Form $L^2$ inner-product between (\ref{L35.1}) and $-e^{2\alpha t}\tilde\Delta\bu_t$ to obtain
\begin{eqnarray}
\label{L37.1}
\frac{\nu}{2}  \frac{d}{dt}\|\tilde\Delta \hat\bu\|^2  
+ e^{2\alpha t}\Big(\|\nabla\bu_t\|^2+\kappa\|\tilde\Delta\bu_t\|^2\Big)
&=&e^{2\alpha t}({\bf f},-\tilde\Delta \bu_t)+e^{2\alpha t}(\bu\cdot\nabla\bu,\tilde\Delta\bu_t)\nonumber\\
&+& \nu\;\alpha \|\tilde\Delta \hat\bu\|^2=I_1 + I_2 +\nu\;\alpha \|\tilde\Delta \hat\bu\|^2.
\end{eqnarray}
Now, integrate (\ref{L37.1}) with respect to time from $0$ to $t$ and then, multiply by  
$2 e^{-2\alpha t}$ to arrive at
\begin{eqnarray}
\label{L37.2}
\nu\|\tilde\Delta\bu\|^2 &+& 2 e^{-2\alpha t} \displaystyle{\int_0^t}e^{2\alpha s}(\|\nabla\bu_t\|^2
+\kappa\|\tilde\Delta\bu_t\|^2)\;ds \leq \nu\;e^{-2\alpha t}\;\|\tilde\Delta\bu_0\|^2\nonumber\\
 &+& 2 e^{-2\alpha t} \int_{0}^t \Big(I_1(s) + I_2(s) \Big)\;ds +2\nu\;\alpha  e^{-2\alpha t} 
 \int_{0}^t e^{2\alpha s}  \|\tilde\Delta \bu(s)\|^2\;ds .
\end{eqnarray}
For $I_2$ on the right hand side of (\ref{L37.1}), rewrite it as 
\begin{align}\label{n3}
I_2= e^{2\alpha t}\;(\bu\cdot\nabla\bu,\tilde\Delta\bu_t)&=\frac{d}{dt} \Big( e^{2\alpha t}\;
(\bu\cdot \nabla\bu,\tilde\Delta\bu)\Big)-2\alpha e^{2\alpha t}\;
(\bu\cdot \nabla\bu,\tilde\Delta\bu) \nonumber\\
&- e^{2\alpha t}\;(\bu_t\cdot \nabla\bu,\tilde\Delta\bu)- e^{2\alpha t}\;(\bu\cdot \nabla\bu_t,\tilde\Delta\bu).
\end{align}
Note that an application of the Ladyzhenskaya's inequality (\ref{LI-d3}) with the Young's inequality shows that
\begin{equation}\label{estimate-LI-1}
e^{2\alpha t}\;
(\bu\cdot \nabla\bu,\tilde\Delta\bu)  \leq  C e^{2 \alpha t}\;\|\bu\|^{1/4}\;\|\nabla \bu\|\;
\|\tilde\Delta\bu\|^{7/4} \leq C(\nu)\; e^{2 \alpha t}\;\|\bu\|^{2}\;\|\nabla \bu\|^8 + 
\frac{\nu}{2} e^{2 \alpha t}  \|\tilde\Delta\bu\|^2.
\end{equation}
From (\ref{L35.1}), we observe using bounds from Lemmas \ref{L43-1} and \ref{L34}  that
\begin{equation}\label{estimate-LI-2}
\|\tilde\Delta\bu\|  \leq  \frac{1}{\nu} \;\Big(\|\bu_t\|+\|\bu\|\;\|\nabla \bu\|+ \|\f\|  
+\kappa \|\tilde \Delta \bu_t\|\Big) \leq C(\nu,\alpha,\lambda_1,M) + \frac{1}{\nu} \kappa \|\tilde\Delta \bu_t\|.
\end{equation}
For the third term on the right hand side of (\ref{n3}), we again employ Ladyzheskaya's 
inequality (\ref{LI-d3}) with estimates from Lemmas  \ref{L43-1}- \ref{L34},(\ref{estimate-LI-2}) and
the Young's inequality
to obtain
\begin{eqnarray}\label{estimate-LI-3}
e^{2\alpha t}\;
(\bu_t\cdot \nabla\bu,\tilde\Delta\bu) 
&\leq& C \; e^{2 \alpha t}\;\|\bu_t\|^{1/4}\;\|\nabla \bu_t\|^{3/4} 
\;\|\tilde\Delta\bu\|^{7/4}\nonumber\\ 
&\leq & C \; e^{2 \alpha t}\;\|\bu_t\|^{1/4}\;\|\nabla \bu_t\|^{3/4} 
\;\Big( C(\nu,\alpha,\lambda_1,M) + \kappa \;\|\tilde\Delta\bu_t\|\Big)^{7/4}\nonumber\\
&\leq&  C(\nu,\alpha,\lambda_1,M)\; e^{2 \alpha t}\;\|\bu_t\|^{1/4}\;\|\nabla \bu_t\|^{3/4} \nonumber\\
&+&
 C(\nu,\alpha,\lambda_1,M)\; e^{\frac{1}{4} \alpha t}\;\|\bu_t\|^{1/4}\; \kappa^{7/8}\|\nabla \bu_t\|^{3/4}\;
 \Big(\sqrt{\kappa}\;\| e^{\alpha t}\tilde\Delta\bu_t\|\Big)^{7/4}\nonumber\\
 &\leq&  C(\nu,\alpha,\lambda_1,M)\; e^{2 \alpha t}\Big(1 + \|\nabla \bu_t\|^2\Big) \nonumber\\
 &+&
  C(\nu,\alpha,\lambda_1,M)\; e^{2 \alpha t}\;\|\bu_t\|^{2}\; \kappa^{4}\; 
  \Big(\kappa\;\|\nabla \bu_t\|^{2}\Big)^3
+\frac{1}{4}  e^{2 \alpha t} \kappa\; \|\tilde\Delta\bu_t\|^2\nonumber\\
&\leq&  C(\nu,\alpha,\lambda_1,M)\; e^{2 \alpha t}\Big(1 + \|\nabla \bu_t\|^2 +\kappa^{4}\;\|\bu_t\|^{2}\Big)  
+\frac{1}{4}  e^{2 \alpha t} \kappa\; \|\tilde\Delta\bu_t\|^2.
\end{eqnarray}
Moreover for the last term  on the right hand side of (\ref{n3}), a use of  following Agmon  inequality
(see, \cite{Foias2001} which is valid for 3D)
\begin{equation}\label{sobolev-inequality}
\|\bu\|_{L^{\infty}} \leq C \|\nabla \bu\|^{1/2}\;\|\tilde{\Delta}\bu\|^{1/2},
\end{equation}
with estimates from Lemmas  \ref{L43-1}- \ref{L34},(\ref{estimate-LI-2}) and
the Young's inequality yields
\begin{eqnarray}\label{estimate-LI-4}
e^{2\alpha t}\;
(\bu\cdot \nabla\bu_t,\tilde\Delta\bu)  &\leq&  C e^{2 \alpha t}\;\|\bu\|_{L^{\infty}}\;\|\nabla \bu_t\|\;
\|\tilde\Delta\bu\| \leq C \; e^{2 \alpha t}\;\|\nabla\bu\|^{1/2}\;\|\tilde{\Delta} \bu\|^{1/2}\;
\|\nabla \bu_t\|\;  \|\tilde\Delta\bu\|\nonumber\\
&\leq& C\; e^{2 \alpha t}\;\|\nabla\bu\|^{1/2}\;\|\nabla \bu_t\|\; \Big( C(\nu,\alpha,\lambda_1,M) + \kappa \;\|\tilde\Delta\bu_t\|\Big)^{3/2}\nonumber\\
&\leq& C\; e^{2 \alpha t}\;\Big( 1+ \|\nabla \bu_t\|^2\Big) +C(\nu,\alpha,\lambda_1,M)\; e^{2 \alpha t}\;
\|\nabla \bu_t\| \;\kappa^{3/2}\|\tilde\Delta\bu_t\|^{3/2}\nonumber\\
&\leq& C\; e^{2 \alpha t}\;\Big( 1+ \|\nabla \bu_t\|^2\Big) +C\;\kappa\; e^{2 \alpha t}\;
\Big(\kappa \|\nabla \bu_t\|^2\Big)^2 + \frac{1}{4}\kappa \|\tilde\Delta\bu_t\|^{2}.
\end{eqnarray}
Substituting (\ref{estimate-LI-3}) and (\ref{estimate-LI-4}) in $I_2$ and integrating with respect to time, use 
{\it a priori} bounds in Lemmas \ref{L43-1}- \ref{L34} to  arrive for the second
term on the right hand side of (\ref{L37.2})  at
\begin{eqnarray}\label{I-2-estimate}
2 e^{-2\alpha t} \int_{0}^t I_2(s)\;ds &\leq& C (\nu,\alpha,\lambda_1,M) + C\;e^{-2\alpha t} \int_{0}^t 
e^{2\alpha s}\Big(1+ (1+\kappa)\|\nabla \bu_t\|^2+ \|\bu_t\|^2 
+ \|\tilde{\Delta}\bu\|^2\Big)\;ds\nonumber\\
&+&\frac{\nu}{4} \|\tilde\Delta\bu(t)\|^2 + e^{-2 \alpha t} \int_0^{t} e^{2\alpha s}\Big(\|\nabla\bu_t\|^2
+\kappa\|\tilde\Delta\bu_t\|^2\Big)\;ds \nonumber\\
&\leq & C(\nu,\alpha,\lambda_1,M)\; +\frac{\nu}{4} \|\tilde\Delta\bu(t)\|^2 
+ e^{-2 \alpha t} \int_0^{t} e^{2\alpha s}\Big(\|\nabla\bu_t\|^2 +\kappa\|\tilde\Delta\bu_t\|^2\Big)\;ds. 
\end{eqnarray}
For $I_1$ term, again rewrite it 
\begin{align}\label{n3-1}
I_1=e^{2\alpha t}\;(\f,\tilde\Delta\bu_t)=\frac{d}{dt}\big(e^{2\alpha t}(\f,\tilde \Delta \bu)\Big)
-2 \alpha e^{2\alpha t}(\f,\tilde \Delta \bu)-e^{2\alpha t}(\f_t,\tilde \Delta \bu).
\end{align}
Now  integrate $I_1$ with respect to time and then multiply by $2 e^{-2\alpha t}$. Then,  a use of 
assumption {\bf (A2)} shows
\begin{eqnarray}\label{n3-2}
2 e^{-2\alpha t} \int_{0}^t I_1(s)\;ds &=& (\f,\tilde \Delta \bu)-e^{-2\alpha t} (\f_0,\tilde \Delta \bu_0)
\nonumber\\
&-& 2 e^{-2\alpha t} \int_{0}^t \alpha e^{2\alpha s}
\Big( 2\alpha (\f,\tilde \Delta \bu) + (\f_t,\tilde \Delta\bu)\Big)\;ds\nonumber\\
&\leq& C(M) + \frac{\nu}{4} \|\tilde\Delta\bu(t)\|^2 + C(\alpha) e^{-2\alpha t} \int_{0}^t \alpha e^{2\alpha s}
\Big(\|\f\|^2+ \|\f_t\|^2\Big)\;ds\nonumber\\
&+& C e^{-2\alpha t} \int_{0}^t \alpha e^{2\alpha s}
\|\tilde{\Delta}\bu\|^2\;ds.
\end{eqnarray}
Substitute (\ref{I-2-estimate}) and (\ref{n3-2}) in (\ref{L37.2})  and use Lemmas 
\ref{L42}, \ref{L43-1}-\ref{L34} with assumption {\bf (A2)} and standard kickback argument 
to arrive at the desired estimate (\ref{estimate-Delta}). To prove (\ref{estimate-Delta-1}),  we note from
(\ref{L35.1}) using Lemmas Lemmas \ref{L43-1}- \ref{L34} with estimate (\ref{LI-d3}) and 
(\ref{estimate-Delta}) that 
\begin{eqnarray*}
\kappa\|\tilde\Delta\bu_t(t)\| &\leq& \|\bu_t\| + \nu \|\tilde\Delta\bu\| + \|\bu\;\cdot \nabla \bu\| + \|\f\|
\nonumber\\
&\leq& \Big(\|\bu_t\| + \nu \|\tilde\Delta\bu\| 
+ C\|\bu\|^{1/4}\;\|\nabla\bu\|\;\|\tilde{\Delta} \bu\|^{3/4} + \|\f\|\Big)
\leq C.
\end{eqnarray*}
 This completes the rest of the proof.\hfill{$\Box$}

The following Lemma \ref{L38} deals with {\it a priori} bounds of the  pressure term. 
\begin{ldf}\label{L38}
Under
assumptions {\rm{({\bf A1})-({\bf A2})}}, there exists
a positive constant $C=C(\nu, \lambda_1,\alpha, M)$ such that for $\displaystyle{0<\alpha< \frac{\nu \lambda_1}{4\left(1+\lambda_1\kappa\right)}}$   and
for all $t>0$, the following estimate holds true:
\begin{eqnarray*}
\| p(t) \|^2_{L^2 /{\rm I\!R}}+ \| p(t)\|^2_{H^1 /{\rm I\!R}}+ e^{-2\alpha t}\int_0^t 
e^{2\alpha s} \| p(s)\|^2_{H^1 /{\rm I\!R}}ds \leq C.
\end{eqnarray*}
\end{ldf}
\noindent
{\it Proof}. A use of the Cauchy-Schwarz inequality with the H\"{o}lder inequality and 
(\ref{LI-d3}) in (\ref{2.2}) yields 
\begin{eqnarray}\label{L38.1}
(p,\nabla\cdot\bphi)\leq C\;\big(\|\bu_t\|+ \kappa\|\nabla \bu_t\| + \|\nabla\bu \|
+\|\nabla\bu\|^2+\|\bf f\|\big)
\|\nabla\bphi\|. 
\end{eqnarray}
 Divide (\ref{L38.1}) by $\|\nabla\bphi\|$ and apply continuous inf-sup condition in (\ref{L38.1}) to obtain
\begin{eqnarray}
\label{L38.2}
\|p\|_{L^2 /{\rm I\!R}}\leq  \frac{|\big(p,\nabla .\bphi\big)|}{\|\nabla\bphi\|}  \leq C\; 
\big( \|\bu_t\|+ \kappa \|\nabla \bu_t\|+\|\nabla\bu \|+\| \nabla \bu \|^{2}+\|\bf f\|\big).
\end{eqnarray}
An application of Lemmas \ref{L42},\ref{L34} and assumption {\bf (A2)} in (\ref{L38.2}) shows
\begin{eqnarray}
\label{L38.3}
\| p(t) \|_{L^2/{\rm I\!R}}\leq C(\nu, \lambda_1, \alpha, M).
\end{eqnarray}
 Use the property of space ${\bf J}_1$ (see \cite{temam} page no 19, remark 1.9) in (\ref{2.3}) to arrive at 
\begin{eqnarray}
\label{L38.4}
(\nabla p,~\bphi)= (\bu_t- \kappa \tilde \Delta\bu_t-\nu\tilde\Delta\bu+\bu.\nabla\bu-\bf f,~\bphi)\,\,\,\,\,\, 
\forall \bphi \in {\bf J}_1.
\end{eqnarray}
A use of the Cauchy-Schwarz inequality with the H\"{o}lder inequality and (\ref{LI-d3}) in (\ref{L38.4}) yields  
\begin{eqnarray}
\label{L38.5}
|(\nabla p,~\bphi)|\leq C(\nu) \left(\|\bu_t\|+ \kappa \|\tilde \Delta\bu_t\|+\|\tilde\Delta\bu\|
+\|\nabla\bu\|\|\tilde\Delta\bu\|^{3/4}
+\|{\bf f}\|\right)\|\bphi\|,
\end{eqnarray}
 and hence,
\begin{eqnarray}
\label{L38.6}
\|\nabla p\|\leq C(\nu)\big(\|\bu_t\|+ \kappa \|\tilde\Delta\bu_t\|+\|\tilde\Delta\bu\|
+ \|\nabla\bu\|\;\|\tilde\Delta \bu\|^{3/4}+\|\bf f\| \big).
\end{eqnarray}
A use of Lemmas \ref{L43-1}, \ref{L34} and \ref{L37} in (\ref{L38.6}) yields
\begin{eqnarray}
\label{L38.7}
\| p(t)\|_{H^1 /{\rm I\!R}} \leq C.
\end{eqnarray} 
Take square of both sides of (\ref{L38.6}). Then, multiply the resulting equation by $e^{2\alpha t}$ 
and integrate from $0$ to $t$ with respect to time to obtain
\begin{align}
\label{L38.8}
\displaystyle{\int_0^t} e^{2\alpha s}\|\nabla p(s)\|^2 \;ds&\leq C(\nu)\bigg(\displaystyle{\int_0^t} e^{2\alpha s}
\big(\|\bu_t(s)\|^2+ \kappa \|\tilde\Delta\bu_t(s)\|^2\big)\;ds
 +\displaystyle{\int_0^t}e^{2\alpha s}\big(\|\tilde\Delta\bu(s) \|^2\nonumber\\
&\qquad+ \|\nabla\bu(s)\|\;\|\tilde\Delta \bu(s)\|^{3/4} \big)\; ds
+\int_0^te^{2\alpha s}\|{\bf f}(s)\|^2 \;ds\bigg).
\end{align}
An application of Lemmas \ref{L43-1}, \ref{L35} and \ref{L37} leads to
\begin{align}
\label{L38.9}
e^{-2\alpha t}\displaystyle{\int_0^t} e^{2\alpha s}\|\nabla p(s)\|^2 ds\leq C.
\end{align}
A use of (\ref{L38.3}), (\ref{L38.7}) and (\ref{L38.9}) would lead to the desired result. 
This concludes the rest of the proof. \hfill{$\Box$}

The main Theorem of this section is stated below without proof as its  proof  follows easily from   
Lemmas \ref{L42},\ref{L43-1}-\ref{L38}. 
\begin{tdf}\label{T31}
Let the assumptions {\rm{({\bf {A1}})}} and {\rm{({\bf {A2}})}} hold. Then, there exists a 
positive constant $C =C(\nu,\alpha, \lambda_1, M)$ such that 
for $\displaystyle{0\leq\alpha< \frac{\nu \lambda_1}{2\big(1+\lambda_1\kappa\big)}}$ the following 
estimates hold true:
\begin{align}
&\| \bu(t)\|_2^2+\| p(s)\|^2_{L^2 /{\rm I\!R}}+ e^{-2\alpha t}
\int_0^te^{2\alpha s}(\|\bu(s)\|^2_2+ \| p(s)\|^2_{H^1 /{\rm I\!R}})\,ds\leq C,\nonumber\\ 
& \|\bu_t(t)\|^2+\kappa \| \bu_t(t)\|_1^2+ \| p(s)\|^2_{H^1 /{\rm I\!R}} +e^{-2\alpha t}
\int_0^t e^{2\alpha s}(\|\bu_{s}(s)\|_1^2+\kappa \| \bu_{s}(s)\|_2^2)ds\leq C. \nonumber
\end{align}
\end{tdf}
\noindent
\begin{remark}
Results in the Theorem \ref{T31} are valid uniformly for all time $t>0$ and even for small $\kappa$ in  
$2D$ and  for  $3D$ with data small. As a result, we can take limit of the equations (\ref{2.2}) 
as $\kappa$ tends to zero which may result in the convergence of the Kelvin-Voigt system to 
the Navier-Stokes system.

Note that an application of Lemmas \ref{L42},\ref{L43}-\ref{L38} instead of Lemma \ref{L43-1} would 
easily provide results of Theorem \ref{T31}, which are valid for both $2D$ and $3D$ without data small, but
 with constant $C$ in the Theorem \ref{T31}   now depending on $1/\kappa.$
 \end{remark}
\begin{remark}
 If $\f \in L^2(0,\infty;\bL^2)$, Theorem \ref{T31} holds uniformly in time with $\alpha=0$.
 When $f(t) =O(e^{-\alpha_0 t}),$ then simple modifications in all Lemmas show exponential decay property 
  which is of order $O(e^{-\alpha_1 t}),$ where $ \alpha_1 =\min(\alpha, \alpha_0)$ in Theorem \ref{T31}.
\end{remark}
\section{ \normalsize \bf The semidiscrete  scheme}
\setcounter{tdf}{0}
\setcounter{ldf}{0}
\setcounter{cdf}{0}
\setcounter{equation}{0}
With $h>0$ as a discretization parameter, let ${\bf H}_h$ and $L_h$, $0<h<1$ be finite dimensional subspaces of 
${\bf H}_0^1 $ and $L^2$, respectively, and be such that, there exist operators 
$i_h$ and $j_h$ satisfying the following approximation properties:  \\
\noindent
({\bf B1}). For each $\bv \in {\bf {J}}_1 \cap {\bf {H}}^2 $
and $ q \in H^1/ {\rm I\!R}$, there are approximations $i_h \bv \in {\bf
{J}}_h $ and $ j_h q \in L_h $ such that
\be
\|\bv-i_h\bv\|+ h \| \nabla (\bv-i_h \bv)\| \le K_0 h^2
\| \bv\|_2, \;\;\;\; \| q - j_h q
 \|_{L^2 /{\rm I\!R}} \le K_0 h \| q\|_{H^1 / {\rm I\!R}}.
\ee
\noindent
For defining the Galerkin approximations, for
$\bv, \bw, \bphi \in \bH_0^1$, set
$
a(\bv, \bphi) = (\nabla \bv, \nabla \bphi)
$
and $b(\bv, \bw,\bphi)$ as in Section 2.
Note that, the operator $b(\cdot, \cdot, \cdot)$ preserves the
antisymmetric properties of the original nonlinear term, i.e.,
$$
b(\bv_h, \bw_h, \bw_h) = 0 \;\;\; \forall \bv_h, \bw_h \in
{\bH}_h.
$$
\noindent
The discrete analogue of the weak formulation (\ref{2.2}) is to find $\bu_h(t) \in {\bf H}_h$ and $p_h(t) \in L_h$ such that $ \bu_h(0)
 = \bu_{0h} $ and for $t>0$,
\begin{align}
\label{4.1}
(\bu_{ht}, {\bphi}_h)+\kappa a ( \bu_{ht}, {\bphi}_h) +\nu a ( \bu_h, {\bphi}_h) &+ b( \bu_h, \bu_h, {\bphi}_h)
-(p_h, \nabla \cdot {\bphi}_h)
=(\bf f, {\bphi}_h) \;\;\; \forall {\bphi}_h \in {\bf H}_h,\nonumber\\
(\nabla \cdot \bu_h, \chi_h) &= 0 \;\;\; \forall \chi_h \in L_h,
\end{align}
\noindent
where $\bu_{0h} \in {\bf H}_h $ is a suitable approximation of $\bu_0\in
{\bf J}_1$ to be defined later.\\
\noindent
We now introduce  $\bJ_h$ as
\be
{\bf J}_h = \{ \bv_h \in {\bf H}_h : (\chi_h, \nabla \cdot \bv_h)
=0 \;\;\; \forall \chi_h \in L_h \}.
\ee
Note that, $\bJ_h$ is not a  subspace of $\bJ_1$. 
Now, the semidiscrete approximation in  $\bJ_h$  is to seek 
 $\bu_h (t) \in {\bf J}_h $ such that $\bu_h(0) = \bu_{0h}\in {\bf J}_h$
and for $t>0$
\ben
\label{4.2}
(\bu_{ht},{\bphi}_h) +\kappa a (\bu_{ht},{\bphi}_h)
+\nu a (\bu_h,{\bphi}_{h})
 =- b( \bu_h, \bu_h, {\bphi}_h)+(\bf f, {\bphi}_h) 
 \;\;\; \forall {\bphi}_h \in {\bf J}_h.
\een
Since $\bJ_h$ is finite dimensional, the equation (\ref{4.2}) leads to a system of 
nonlinear ordinary differential equations. Therefore, an application of Picard's theorem ensures 
existence of a unique solution $\bu_h$ for $(0,t_h^{*})$ for some $t_h^{*}>0$. For global existence, 
we need to use continuation 
argument  provided the discrete solution is bounded  for all $t>0.$
Following the argument in the proof of Lemma \ref{L42}, it is easy to prove the following estimate:
for  $\displaystyle{0<\alpha<\frac{\nu\lambda_1}{4\left(1+\kappa\lambda_1 \right)}}$ and  for all $t > 0$
\begin{align}
 \label{u-h-estimate}
\|\bu_h(t)\|^2 &+\kappa\|\nabla \bu_h(t)\|^2+ \beta e^{-2\alpha t}
\int_0^t{e^{2\alpha s}\|\nabla \bu_h(s)\|^2}\,ds\nonumber \\
&\le  e^{-2\alpha t}(\|\bu_{0h}\|^2+\kappa \|\nabla\bu_{0h}\|^2)+ 
\left(\frac{1-e^{-2\alpha t}}{2\nu\lambda_1\alpha}\right)
\|{\bf f}\|^2_{L^{\infty}(\bL^2)},
\end{align}
where $\beta= \nu-2\alpha (\kappa+\lambda_1^{-1})>\nu/2>0.$ 
This complete the proof of existence and uniqueness of a global discrete solution for all $t>0$.

As a consequence of (\ref{u-h-estimate}), the following result on existence of a discrete 
global attractor is derived.
\begin{ldf}\label{discrete-attractor}
There exists a bounded absorbing set
$$
B_{\brho_0}(0)=\{\bu_h \in \bJ_h : \Big(\|\bu_h\|^2 +\kappa \|\nabla \bu_h\|^2\Big)^{1/2} \leq \brho_0 \}$$ 
with $\brho_0$ given by 
 $$\brho_0^2=\left(\frac{1}{\alpha \nu\lambda_1 }\right)\|\f\|^2_{L^{\infty}(\bL^2)}.
 $$
 Further, the problem (\ref{4.2}) has a global attractor ${\mathcal{A}}_h \subset \bJ_h,$ which attracts bounded
 sets in $\bJ_h.$ 
 \end{ldf}
 \noindent
 {\it Proof}. To prove the first part, we need to show an existence of $\brho_1 >0$ such that 
 for any $\bu_{0h}\in \bJ_h,$ there exists a time $t^*:=t^*((\|\bu_{0h}\|^2 +\kappa 
 \|\nabla \bu_{0h}\|^2)^{1/2})$ 
 such that for $t \geq t^*$ the  discrete solution $\bu_h(t)$ of (\ref{4.2}) satisfies $ \bu_h(t) \in 
 B_{\brho_1}.$ For any ball $B_{\brho_1}(0),\;\;\brho_1 > \brho_0/2$  with
 the initial condition $\bu_{0h}\in B_{\brho_1}(0),$ 
 it follows from (\ref{u-h-estimate}) that
 \begin{eqnarray}
  \label{u-h-1}
 \Big(\|\bu_h(t)\|^2 +\kappa\|\nabla \bu_h(t)\|^2\Big)^{1/2}
 &\le&  e^{-2\alpha t} \brho_1^2 + 
 \frac{1}{2} \brho_0^2 \left(1-e^{-2\alpha t}\right) \\
 &=& e^{-2\alpha t} \left(\brho_1^2- \frac{1}{2} \brho_0^2\right) + 
  \frac{1}{2} \brho_0^2. \nonumber
 \end{eqnarray}
 To complete the proof, we claim that
 $$e^{-2\alpha t} \left(\brho_1^2- \frac{1}{2} \brho_0^2\right) \leq 
   \frac{1}{2} \brho_0^2.
$$
This can be achieved if 
$$t  \geq \frac{1}{\alpha} \log \Big( \frac{2 \brho_1^2 - \brho_0^2}{\brho_0^2}\Big)=:t^* >0,$$
that is, for $t\geq t^*,$ $B_{\brho_1}(0) \subset B_{\brho_0}(0).$
Note that for  $\brho_1 \leq \brho_0/2,$ it is trivially satisfied for all $t>0.$
Hence, $B_{\brho_0}(0)$ is an absorbing ball and it further follows that the problem (\ref{4.2}) has 
a discrete global  attractor  ${\mathcal{A}}_h \subset \bJ_h,$ which attracts bounded
 sets in $\bJ_h.$  This completes the rest of the proof. \hfill{$\Box$}
 
Define the quotient space $L_h/N_h$, where
\be
N_h = \{ q_h \in L_h : (q_h, \nab \cdot {\bphi}_h ) = 0, \forall
{\bphi}_h \in {\bf H}_h\}
\ee
with its norm given by
\be
\| q_h\|_{L^2/N_h} = {\dps{\inf_{\chi_h \in N_h} }} \| q_h + \chi_h\|.
\ee
Furthermore, assume that the pair $({\bf H}_h,{L_h/N_h})$ satisfies the following  uniform inf-sup condition:\\
\noindent
({\bf B2}). For every $q_h \in L_h$, there exist a non-trivial
function $\bphi_h \in \bf H_h$ and a positive constant $K_1$, independent of $h$, such that
\be
|(q_h, \nab \cdot {\bphi}_h)| \ge K_1\|\nab {\bphi}_h \| \| q_h\|_{
L^2/N_h}.
\ee
As a consequence of conditions ({\bf B1})-({\bf B2}), we have 
the following properties of 
the $L^2$ projection $P_h:\bL^2\rightarrow \bJ_h$. 
For $\bphi \in \bJ_1$, we note that, (see  \cite{GR}, \cite{HR82}), 
\ben
\label{4.4}
\|\bphi- P_h \bphi\|+ h \|\nabla P_h \bphi\| \leq C h
\|\nabla \bphi\|,
\een
and for $\bphi \in \bJ_1 \cap \bH^2$
\ben
\label{4.5}
\|\bphi- P_h \bphi\|+ h \|\nabla (\bphi-P_h \bphi)\| \leq C h^2
\|\tilde \Delta \bphi\|.
\een
We may define the discrete  operator $\Delta_h: \bH_h
\rightarrow \bH_h$ through the bilinear form $a (\cdot, \cdot)$ as
\begin{eqnarray}
\label{4.6}
a(\bv_h, {\bphi}_h) = (- \Delta_h \bv_h, \bphi) \;\;\;\;
\forall \bv_h, {\bphi}_h \in \bH_h.
\end{eqnarray}
Set the discrete analogue of the Stokes operator $\tl{\Delta} =
P \Delta $ as
 $\tl{\Delta}_h = P_h \Delta_h $.
Examples of subspaces $\bf H_h$ and $L_h$ satisfying assumptions $(\rm{\bf {B1}})$ and $(\rm{\bf {B2}})$ can be 
found in \cite{BP} and \cite{HR82}. \\
\noindent
Next in the following Lemma, {\it a priori} bounds for the discrete 
solution $\bu_h$ of (\ref{4.2}), which will be helpful in establishing the error estimates, are stated. 
The proof can be obtained following the similar steps 
as in the proofs of Lemma 
\ref{L42}-\ref{L35}. 	
\begin{ldf} 
\label{L41}
For all $t>0$, the semi-discrete Galerkin approximation
$\bu_h$ for the velocity satisfies 
\begin{eqnarray}
\|\bu_{h}(t)\|_1^2+\kappa\|\tilde\Delta_h\bu_{h}(t)\|^2+\|\tilde\Delta_h\bu_{h}(t)\|^2+e^{-2\alpha t}
\int_0^t e^{2\alpha s}
(\|\nab\bu_{h}\|^2+\|\tilde\Delta_h\bu_{h}\|^2+\|\nab\bu_{ht}\|^2)\,ds \le C.\nonumber
\end{eqnarray}
\end{ldf}
\section{ \normalsize \bf Error estimates for the velocity}
\setcounter{tdf}{0}
\setcounter{ldf}{0}
\setcounter{cdf}{0}
\setcounter{equation}{0}
In this section, we analyze the error occurred due to the Galerkin approximation for the velocity term. 

Since $\bJ_h$ is not a subspace of $\bJ_1$, the weak solution $\bu$ satisfies
\begin{eqnarray}
\label{E62}
(\bu_t,{\bphi}_h)+\kappa a(\bu_{t},{\bphi}_h)+\nu a(\bu,{\bphi}_h)=-b(\bu,\bu,{\bphi}_h)+(p,\nabla \cdot {\bphi}_h)+({\bf f},{\bphi}_h) 
\,\,\,\, \forall {\bphi}_h \in \bJ_h.  
\end{eqnarray}
Set ${\bf e}=\bu-\bu_h$. Then, from (\ref{E62}) and (\ref{4.2}), we obtain
\begin{eqnarray}\label{E63}
({\bf e}_{t},{\bphi}_h)+\kappa a({\bf e}_{t},{\bphi}_h)+\nu a({\bf e},{\bphi}_h)= {\bf \Lambda}({\bphi}_h)
+(p,\nabla \cdot {\bphi}_h),
\end{eqnarray}
\noindent
where ${\bf \Lambda}({\bphi}_h)=-b(\bu,\bu,{\bphi}_h)+b(\bu_h,\bu_h,{\bphi}_h)$.
Below, we derive an optimal error estimate of $\|\nabla\e(t)\|$, for $t>0$.
\begin{ldf}
\label{T61}
 Let assumptions {\rm (\bf{A1})-(\bf{A2})} and {\rm (\bf{B1})-(\bf{B2})} be satisfied.  With  
 $\bu_{0h}= P_h \bu_0,$ then, there exists a positive constant $C$ depending on $\lambda_1$, $\nu$, $\alpha$ 
 and $M$, such that, for fixed $T > 0$ with  $t \in (0, T)$ and for $\displaystyle{0\leq\alpha < \frac{\nu \lambda_1}{4\big(1+\lambda_1\kappa\big)}}$, the following estimate holds true :
\begin{align}
\|(\bu-\bu_h)(t)\|^2 + \kappa \|\nabla(\bu -\bu_h)(t)\|^2 \leq C h^2 e^{CT}.\nonumber
\end{align}
\end{ldf}
\noindent
{\it Proof.} On multiplying(\ref{E63}) by $e^{\alpha t} $ with 
${\bphi}_h=P_h\hat{\bf e}=\hat\e+(P_h \hat \bu-\hat \bu),$ it follows that
\begin{align}
\label{E64}
(e^{\alpha t}\e_t,\hat{\e}) &+\kappa a(e^{\alpha t}\e_t,\hat{\e}) + \nu a(\hat{\e},\hat{\e}) =e^{\alpha t}\Lambda(P_h\hat{\bf e}) + (\hat{p},\nabla\cdot P_h\hat {\bf e})\nonumber\\
&+(e^{\alpha t}\e_t,\hat{\bu}-P_h\hat{\bu}) + \kappa\, a(e^{\alpha t}\e_t,\hat{\bf u}-P_h\hat{\bf u})+\nu a(\hat{\e},\hat{\bf u}-P_h\hat{\bf u}).
\end{align}
Note that
\begin{align}\label{E64*}
(e^{\alpha t}{\bf e}_t,\hat{\bf e})+\kappa\,a( e^{\alpha t}{\bf \e}_t,\hat{\bf e})=\displaystyle\frac{1}{2}\frac{d}{dt}(\|\hat{\bf e}\|^2+\kappa \|\nabla\hat{\bf e}\|^2) -\alpha (\|\hat{\bf e}\|^2+\kappa \|\nabla\hat{\bf e}\|^2),
\end{align}
\noindent
and using $L^2$-projection $P_h,$ we find that
\begin{align}\label{E65*}
 (e^{\alpha t}\e_t,\hat{\bu}-P_h\hat{\bu})&=(e^{\alpha t}({\e}_t-P_h \e_t) ,\hat{\bu}-P_h\hat{\bu}))
 -\alpha(e^{\alpha t}({\e}-P_h \e_t),\hat{\bu}-P_h\hat{\bu})\nonumber\\
 &=\frac{1}{2}\frac{d}{dt} \| \hat{\bu}-P_h\hat{\bu}\|^2 -\alpha \|\hat{\bu}-P_h\hat{\bu})\|^2.
\end{align}
 A use of (\ref{2.1*}) with (\ref{E64*}) and (\ref{E65*}) in (\ref{E64}) yields  
\begin{align}\label{E65}
\displaystyle\frac{d}{dt}(\|\hat{\bf e}\|^2  &+  \kappa \|\nabla\hat{\bf e}\|^2)
+ \left( 2\nu  -2\alpha(\kappa+{\lambda_1}^{-1} 
)\right) \|\nabla\hat{\bf e}\|^2\leq 2e^{\alpha t} \Lambda(P_h\hat{\bf e}) 
+ 2(\hat{p},\nabla \cdot P_h\hat{\bf e})\nonumber\\
& + \frac{d}{dt}\bigg(\|\hat{\bu}-P_h\hat{\bu}\|^2 + 2 \kappa\,a(\hat{\bf e},\hat{\bu}-P_h\hat{\bu})\bigg)
- 2 \kappa\,a(\hat{\bf e}, e^{\alpha t}({\bu}_t-P_h\hat{\bu}_t))\nonumber\\
&-2\alpha \bigg( \|\hat{\bu}-P_h\hat{\bu}\|^2 +\kappa\,a(\hat{\bf e},\hat{\bu}-P_h\hat{\bu})\bigg)
 +2\nu a(\hat{\bf e},\hat{\bu}-P_h\hat{\bu}).
\end{align}
\noindent
For the last three terms on the right hand side of (\ref{E65}), apply the Cauchy-Schwarz inequality 
with  Poincar\'e inequality and  Young inequality to bound it as
\begin{align}\label{E66}
|2\alpha \big(\|\hat{\bu}-P_h\hat{\bu}\|^2&+\kappa\,a(\hat{\bf e},\hat{\bu}-P_h\hat{\bu})\big)
+2\nu a(\hat{\bf e},\hat{\bu}-P_h\hat{\bu}) + 2 \kappa\,a(\hat{\bf e}, e^{\alpha t}({\bu}_t-
P_h {\bu}_t)|\nonumber\\
&\leq C(\alpha,\lambda_1,\nu,\epsilon) \Big(\|\nabla (\hat{\bu}-P_h\hat{\bu})\|^2 + 
\kappa^2\|e^{\alpha t}\nabla (\bu_t-P_h \bu_t)\|^2
 + \frac{\epsilon}{2} \|\nabla \hat{\bf e}\|^2.
\end{align}
For the second term on the right-hand side of (\ref{E65}), a use of approximation property$({\bf B}_1)$ 
with discrete in compressibility condition and $\bH_0^1 $- stability of the $\bL^2$- projection $P_h$ shows 
\begin{align}
\label{E67}
2|(\hat{p}, \nabla\cdot P_h\hat{\bf e})| &= |(\hat{p}-j_h\hat{p}, \nabla\cdot P_h\hat{\bf e})|
\leq C\|\hat{p}-j_h \hat{p}\| \;\| \nabla P_h\hat{\bf e}\| \nonumber\\
&\leq  C(\epsilon)h^2\;\|\nabla\hat{p}\|^2+\frac{\epsilon}{2} \;\| \nabla \hat{\bf e}\|^2.
\end{align}
To estimate the first term on the right-hand side of (\ref{E65}), use anti-symmetric property (\ref{ASP}) 
of the trilinear form $b(\cdot,\cdot,\dot)$ and the property of $P_h$ to obtain
\begin{align}\label{lambda-estimate}
2e^{\alpha t} \Lambda (P_h\hat{\bf e})=-2e^{-\alpha t}\bigg( b(\hat{\bf e}, \hat{\e},\hat{\bu}-P_h\hat{\bu}) + b(\hat{\bf e}, \hat{\bu},P_h\hat{\bf e})+ b(\hat{\bu},\hat{\bf e},P_h\hat{\bf e})\bigg).
\end{align}
Then,
using the generalized H\"{o}lder inequality, the Agmon inequality (\ref{sobolev-inequality}),
the Young inequality, the Sobolev embedding theorem, (\ref{2.1a}) and (\ref{4.4}), 
we arrive at
\begin{align}\label{E68}
2e^{\alpha t} |\Lambda (P_h\hat{\bf e})| &
\leq 2 e^{-\alpha t}
\big(\| \hat \bu\|_{L^\infty}\|\nabla\hat{\bf e}\|\|P_h \hat {\bf e}\|
+\| \nabla\hat {\bf e}\|\|\tilde\Delta\hat \bu\|\|P_h\hat {\bf e}\| + \|\nabla\hat \e\|\;
\|\nabla\hat\e\|\|\nabla(\hat \bu-P_h \hat{\bu})\| \big)\nonumber\\
&\leq 2 e^{-\alpha t}\Big(\big(\displaystyle{\|\nabla \hat \bu\|^\frac{1}{2}}\displaystyle{\|\tilde\Delta \hat\bu\|}^{\frac{1}{2}}+ \|\tilde \Delta\hat \bu\|\big)\|\hat {\bf e}\|\|\nabla \hat\e\| + (\|\nabla\hat \bu\|+\|\nabla \hat  \bu_h\|) \|\nabla\hat\e\|\|\nabla(\hat \bu-P_h \hat\bu)\|\Big)  \nonumber\\
&\leq C(\epsilon) e^{-2\alpha t}\Big(\big( \|\nabla\hat \bu\|\|\tilde\Delta \hat\bu\|
+\|\tilde\Delta\hat \bu\|^2 \big)\|\hat{\e}\|^2+ \|\nabla (\hat\bu-P_h\hat\bu)\|^2\Big)
+\frac{\epsilon}{2}\|\nabla \hat \e\|^2. 
\end{align}


\noindent
Integrating (\ref{E65}) with respect to time from $0$ to $t$, use bounds (\ref{E66}), \ref{E67} and 
(\ref{E68}) with $\epsilon=\frac{2 \nu}{3}$, to arrive at
\begin{align}\label{E611}
\|\hat{\bf e}(t)\|^2 &+ \kappa \|\nabla\hat{\bf e}(t)\|^2 +\beta \int_0^t\|\nabla\hat{\bf e}\|^2ds \leq C(\|{\bf e}(0)\|^2 + \|\nabla {\bf e}(0)\|^2)\nonumber\\
  &+C(\alpha,\nu,\lambda_1,M)\bigg(\|\nabla(\hat{\bu}-P_h\hat{\bu})\|^2
 +\int_0^t(\|\nabla(\hat{\bu}-P_h\hat{\bu})\|^2 + \kappa^2\|\nabla(\hat{\bu}_t-P_h\hat{\bu}_t)\|^2\nonumber\\
 &~~~~+\|\nabla\hat{p}\|^2 )ds\bigg)+ C \int_0^t\big(\|\nabla\bu\|\|\tilde{\Delta}\bu\|+\|\tilde\Delta \bu\|^2\big)
 \|\hat{\bf e}\|^2\; ds. 
\end{align}
\noindent
A use of  (\ref{4.5}) and $({\bf B1})$ in (\ref{E611}) yields
\begin{align}
\|\hat{\bf e}(t)\|^2+\kappa ||\nabla\hat{\bf e}(t)\|^2 +\beta \int_0^t\|\nabla\hat{\bf e}\|^2ds&\leq C h^2 \bigg(\|\bu_0\|_2^2+\|\hat{\bu}\|_2^2 + \int_0^t(\|\hat{\bu}\|_2^2 +\|\hat{\bu}_t\|_2^2 +\|\hat{p}(t)\|_{H^1 /{\rm I\!R}}^2) ds\bigg)\nonumber\\
&+C \int_0^t\big(\|\nabla\bu\|\|\tilde{\Delta}\bu\| +\|\tilde\Delta \bu\|^2\big)\big(\|\hat{\bf e}\|^2+ \kappa \|\nabla\hat{\bf e}\|^2\big)ds.\nonumber
\end{align}
\noindent
From the {\it a priori } bounds of $\bu$, $\bu_t$ and $p$ in Theorem \ref{T41}, we arrive  using 
 the Gronwall lemma at
\begin{align}
\|\hat{\e}(t)\|^2+\kappa \|\nabla\hat{\e}(t)\|^2 +\beta \int_0^t\|\nabla\hat{\bf e}\|^2\;ds\leq 
C(\nu,\alpha,\lambda_1,M) h^2 \mbox{exp}\bigg(\int_0^t(\|\tilde{\Delta} \bu\|^2 +\|\nabla\bu\| \|\tilde{\Delta}\bu \| )\;ds \bigg).\nonumber
\end{align}
A use of {\it a priori} bounds given in Lemma \ref{L43-1} yields
\begin{equation}\label{u-integral}
\displaystyle{\int_0^t} \big(\|\nabla\bu\|\;\|\tilde{\Delta}\bu \|+ \|\tilde\Delta \bu\|^2\big) ds  \leq Ct, 
\end{equation}
and hence, we find that 
\begin{eqnarray}
\|(\bu-\bu_h)(t)\|^2 + \kappa \|\nabla(\bu -\bu_h)(t)\|^2 \leq C h^2 e^{Ct}.\nonumber
\end{eqnarray}
This concludes the proof.\hfill {$\Box$}

Observe that the Lemma \ref{T61} provides a suboptimal error estimates for the velocity in $L^{\infty}(\bL^2)$-norm. 
Therefore, in the remaining part of this section, we derive an optimal error estimate for the velocity in $L^{\infty}(\bL^2)$-norm.

Introduce an intermediate solution $\bv_h$ which is a finite element Galerkin 
approximation to a linearized Kelvin-Voigt equation, that is ,$\bv_h$ satisfies   
\ben
\label{E614}
 (\bv_{ht}, {\bphi}_h)+\kappa\,a(\bv_{ht},{\bphi}_h)+\nu a (\bv_h, {\bphi}_h)
=(\bf f,{\bphi}_h)-b(\bu,\bu,{\bphi}_h)
\;\;\;\; \forall {\bphi}_h \in {\bf J}_h,
\een
with $\bv_h(0)=P_h \bu_0.$\\
\noindent
Now, we split ${\bf {e}} $ as
\be
{\bf  {e}} := \bu - \bu_h = (\bu - \bv_h) + ( \bv_h - \bu_h)
= \bxi + \bta.
\ee
Note that $\bxi$ is the error committed by approximating  a linearized Kelvin-Voigt equation (\ref{E614}) and $\bta$ represents  the error
due to the non-linearity in the equation.
Now, subtract (\ref{E614}) from (\ref{E62}) to write an  equation in $\bxi$  as
\begin{eqnarray}\label{E615}
(\bxi_t,{\bphi}_h)+\kappa\,a({\bxi}_t,{\bphi}_h)+\nu a(\bxi,{\bphi}_h)=(p,\nabla \cdot {\bphi}_h)\;\;\;\;\forall {\bphi}_h\in {\bf {J}}_h.
\end{eqnarray}
For deriving optimal error estimates of $\bxi$ in $L^{\infty}(L^2)$ and
$L^{\infty}(H^1)$-norms, we introduce, as in \cite{ANS}, the following
Sobolev-Stokes's  projection  $V_h\bu:[0,\infty)\rightarrow J_h$
satisfying
\begin{eqnarray}\label{E616}
\kappa  a(\bu_t-V_h\bu_t,{\bphi}_h)+ \nu a(\bu-V_h\bu,{\bphi}_h)=(p,\nabla \cdot {\bphi}_h)\; \; \;
\forall {\bphi}_h\in \bJ_h,
\end{eqnarray}
where $V_h\bu(0)=P_h \bu_0.$ In other words, given $(\bu,p),$ find $V_h\bu:[0,\infty)\rightarrow J_h$ satisfying (\ref{E616}). 
Since ${\bf J}_h$ is finite dimensional, for a given $\bu$ the problem (\ref{E616}) leads to a linear 
system of ODEs. Then, an application of  Picard's theorem with continuation argument ensures 
existence of a unique solution in $[0,\infty).$
\noindent
With $V_h\bu $ defined as above, we now split $\bxi$ as
$$
\bxi:=(\bu-V_h\bu)+(V_h\bu-\bv_h)=:\bzeta+\brho.
$$
To obtain estimates for ${\bxi}$, first of all, we state estimates of $ \bzeta $ 
in Lemmas \ref{L61} and \ref {L62*}. Then, we proceed to estimate $\|\brho\|$ and $\|\nabla \brho\|$ 
 in Lemma \ref{L65*}. Combining these results, we obtain estimates for $\bxi$ in $L^{\infty}(\bL^2)$ and 
 $L^{\infty}(\bH_0^1)$-norms in Lemma \ref{L63}. Finally, 
 we derive an estimate for $\bta$ to complete the proof of our main Theorem \ref{T41}.

Below, we briefly state the proofs of the above lemmas. The proofs are along  similar lines as in the proofs of 
Lemmas 5.2-5.7 in \cite{ANS}. The difference occur only in applying {\it a priori} estimates 
as they do not decay exponentially in time. Therefore, in the following proofs, we briefly indicate 
the differences.
\begin{ldf}\label{L61}
Assume that {\rm{({\bf A1})-({\bf A2})}} and {\rm{({\bf B1})-({\bf B2})}} are satisfied. Then, there exists a positive 
constant $C=C(\nu,\lambda_1, \alpha, M)$ such that for 
$\displaystyle{0<\alpha< \frac{\nu\lambda_1}{4\big(1+\kappa\lambda_1\big )}}$, 
the following estimate holds true:
\begin{eqnarray*}
  \kappa \|\nabla \bzeta(t)\|^2 +e^{-2\alpha t}\int_0^t e^{2\alpha s}\|\nabla\bzeta(s)\|^2 ds\leq C h^2 .
\end{eqnarray*}
\end{ldf}
\noindent{\it Proof.} We first multiply (\ref{E616}) by $e^{\alpha t}$ with $\bzeta =\bu-V_h\bu$ and 
then choose ${\bphi}_h=P_h\hat{\bzeta}=\hat\bzeta-(\hat \bu-P_h \hat\bu)$ to arrive at 
\begin{align}\label{E619}
 \kappa \frac{d}{dt}\|\nabla \hat{\bzeta}\|^2 +2(\nu-&\kappa \alpha)\|\nabla \hat{\bzeta}\|^2 = 2\kappa \frac{d}{dt}
a ( \hat\bzeta, \hat\bu-P_h \hat{\bu})-2\kappa\,a (\hat\bzeta, \frac{d}{dt}(\hat\bu-P_h \hat{\bu}))\nonumber\\
&+ 2(\nu-\kappa \alpha)\, a(\hat{\bzeta},
\hat{\bu}-P_h\hat{\bu})
+ 2(\hat p,
\nabla \cdot P_h\hat{\bzeta}).
\end{align}
Integrating  (\ref{E619}) with respect to time from $0$ to $t,$ a use of (\ref{4.4}) along with the 
Youngs inequality yields
\begin{align}
\label{E621}
\kappa\|\nabla \hat{\bzeta}\|^2+(\nu-&\kappa \alpha)\int_0^t\|\nabla \hat\bzeta\|^2 \;ds \leq 
C(\nu,\alpha) \bigg(\|\nabla(\bu_0-P_h \bu_0)\|^2+e^{2\alpha t}\|\nabla(\bu-P_h \bu)\|^2\nonumber\\
&+\int_0^t e^{2\alpha s}\big(\|\nabla (\bu_t-P_h \bu_t)\|^2 
+ \|\nabla (\bu-P_h \bu)\|^2+
  \|\nabla p\|^2 \big)\;ds\bigg).
\end{align}
Now, use  (\ref{4.5}) and {\rm{({\bf B1})}} in (\ref{E621}) to obtain
\begin{align}\label{ne2}
\kappa\|\nabla \hat{\bzeta}\|^2+(\nu-\kappa \alpha)&\int_0^t\|\nabla \hat\bzeta\|^2 ds \leq 
C(\nu, \alpha) h^2\bigg( \|\tilde {\Delta}\bu_0\|^2+e^{2\alpha t}\|\tilde {\Delta}\bu\|^2
+\int_0^t e^{ 2\alpha s}\|\nabla p\|^2 \;ds\nonumber\\
&+\int_0^t e^{2\alpha s}( \|\tilde{\Delta}\bu_t\| ^2 
+ \|\tilde{\Delta} \bu \|^2 )\;ds\bigg).
\end{align}
From {\it a priori} bounds for $\bu$ and $p$ derived in Lemmas \ref{L43}, \ref{L37} and \ref{L38}, 
we arrive at the desired result. This completes the rest of the proof.
\hfill{$\Box$}


\noindent
Below, we state a lemma without proof. The proof can be obtained in a similar fashion as in \cite{ANS} and 
applying now {\it a priori} estimates derived in Theorem \ref{T31}.   
\begin{ldf}\label{L62*}
Under the assumptions {\rm{({\bf A1})-({\bf A2})}} and {\rm{({\bf B1})-({\bf B2})}}, there exists a positive constant 
$C=C(\nu,\lambda_1, \alpha,M)$ such that for 
$\displaystyle{0<\alpha< \frac{\nu\lambda_1}{4\big(1+\kappa\lambda_1\big)}}$, the following estimate holds true for $t >0$:
\begin{eqnarray}
\kappa \|\bzeta(t) \|^2+e^{-2\alpha t}\int_0^t e^{2\alpha s}\left(\|\bzeta(s)\|^2+ \kappa\|\bzeta_t(s)\|^2 
+ \kappa h^2\|\nabla \bzeta_t(s)\|^2\right)
ds\leq C h^4.\nonumber
\end{eqnarray}
\end{ldf}
\noindent
In the following Lemma,  estimates of $\brho$  are derived.
\begin{ldf}\label{L65*}
Under the assumptions {\rm{({\bf A1})-({\bf A2})}} and {\rm{({\bf B1})-({\bf B2})}}, there exists a positive constant 
$C=C(\nu,\lambda_1, \alpha,M)$ such that for $\displaystyle{0<\alpha< \frac{\nu}{4\big(1+\kappa\lambda_1\big)}}$, 
the following estimate holds true:
\begin{eqnarray}
 \kappa (\|\brho\|^2+\kappa \|\nabla \brho\|^2)+2\kappa \beta e^{-2\alpha t}
 \int_0^t e^{2\alpha s} \|\brho(s)\|^2 \;ds \leq C(\nu,\lambda_1,\alpha, M) h^4.\nonumber
 \end{eqnarray}
\end{ldf}
\noindent{\it Proof.} Subtract (\ref{E616}) from (\ref{E615}) and substitute $\bphi_h$ by $e^{\alpha t}\hat\brho$ to obtain
\begin{eqnarray}
\label{E647*}
(e^{\alpha t}\brho_t,\hat \brho)+\kappa\,a(e^{\alpha t}\brho_t, \hat\brho)+ \nu\|\nabla\hat\brho\|^2 = -(e^{\alpha t}\bzeta_t, \hat\brho) \;\;\;\;\; \forall {\bphi}_h\in \bJ_h.
\end{eqnarray}
 Apply the Cauchy-Schwarz inequality, (\ref{2.1*}) with the Young inequality in (\ref{E647*}) 
 and integrate with respect to time from 0 to $t$ 
to arrive at  
\begin{align}\label{E6941*}
\|\hat \brho\|^2+\kappa\|\nabla\hat \brho\|^2+2\beta\displaystyle{\int_0^t}\|\nabla \hat \brho\|^2 ds\leq 
C(\alpha,\lambda_1)\int_0^t\| e^{\alpha s}\bzeta_t(s)\|^2ds.
\end{align}
The desired result follows after a use of Lemma \ref{L62*} in (\ref{E6941*}).\hfill{$\Box$}\\
\noindent
We now derive an estimate of $\bxi$ in $L^{\infty}(\bL^2)$ and $L^{\infty}(\bH^1_0)$-norms.
\begin{ldf}\label{L63}
Let the assumptions {\rm{({\bf A1})-({\bf A2})}} and {\rm{({\bf B1})-({\bf B2})}} be satisfied. Then, 
there exists a positive constant $C=C(\nu,\lambda_1, \alpha,M)$ such that for 
$\displaystyle{0<\alpha< \frac{\nu\lambda_1}{4\big(1+\kappa \lambda_1\big)}}$, the following estimate holds:
\begin{eqnarray}
 \kappa \|\bxi(t)\|^2+ \kappa \|\nabla\bxi(t)\|^2+e^{-2\alpha t}\int_0^t e^{2\alpha s} \|\bxi(s)\|^2\; ds 
\leq C(\nu,\lambda_1,\alpha,M) h^4.\nonumber
 \end{eqnarray}
\end{ldf}
\noindent{\it Proof.} A use of the triangle inequality along with Lemmas \ref{L61}-\ref{L65*} leads to the desired result.\hfill{$\Box$}
\begin{ldf}
\label{L64}
Let the assumptions {\rm{({\bf A1})-({\bf A2})}} and {\rm{({\bf B1})-({\bf B2})}} hold true. Let 
$\bu_h(t) \in {\bf J}_h$ be 
a solution of (\ref{4.2}) with initial condition $\bu_h(0)=P_h\bu_0$, where $\bu_0 \in{\bf J}_1 $. 
Then there exist a constant C such that for $0 < T < \infty$ with $t \in (0,T]$
\ben
e^{-2\alpha t}\int_0^t e^{2\alpha t}\|{\bf e}\|^2 \leq C e^{C T} h^4.\nonumber
\een
\end{ldf} 
\noindent{\it Proof.} In view of Lemma \ref{L63}, we only need to prove the estimate for $\bta$. From (\ref{E614}) and (\ref{4.2}), 
the equation in $\bta$ becomes 
\begin{eqnarray} \label{E653}
(\bta_t,{\bphi}_h)+ \kappa\,a(\bta_t,{\bphi}_h) +\nu
a(\bta,{\bphi}_h)= \Lambda_h(\bphi_h),
\,\,\,\forall {\bphi}_h\in\bJ_h,
\end{eqnarray}
 where 
 \ben \label{Lambda}
 \Lambda_h(\bphi_h)= b(\bu_h,\bu_h,\bphi_h)-b(\bu,\bu,\bphi_h)= - b({\bf e},\bu_h,\bphi_h)-b(\bu,{\bf e},\bphi_h).
 \een
Substitute ${\bphi}_h=e^{2\alpha t}({\tilde\Delta}_h^{-1}\bta)$ in (\ref{E653}) to obtain
\ben
\label{L64.1}
\frac{1}{2}\frac{d}{dt}(\|\hat{\bta}\|_{-1}^2 +\kappa \| \hat{\bta}\|^2)-\alpha\|\hat{\bta}\|_{-1}^2+(\nu-\kappa\alpha)\|\hat{\bta}\|^2=e^{\alpha t}\Lambda_h(\hat{\bta}).
\een
We recall that $\|\bw_h\|_{-1} := \|(-{\tilde\Delta}_h)^{-1/2}\bw_h\| $ for $w_h \in {\bf J}_h$. 
Again for $\bv \in {\bf J}_1$ and $\bphi,\,\bxi \in {\bf J}_h$
\ben
\label{L64.1*}
|b(\bv,\bphi,\bxi)| \leq C \|\bv\|^{1/2}\|\nabla\bv\|^{1/2}\|\bphi\|\|\nabla\bxi\|^{1/2}\|{\tilde\Delta}_h\bxi\|^{1/2}.
\een
For $\bv ,\,\bphi,\,\bxi \in {\bf J}_h$
\ben
\label{L64.1**}
|b(\bv,\bphi,\bxi)| \leq \|\bv\|\|\nabla\bphi\|^{1/2}\|{\tilde\Delta}_h\bphi\|^{1/2}(\|\bxi\|^{1/2}\|\nabla\bxi\|^{1/2}+\|\nabla\bxi\|).
\een
Now, a use of ${\bf e}=\bxi+\bta$, along with (\ref{L64.1*}) and (\ref{L64.1**}) leads to 
\ben
\label{L64.1***}
|e^{\alpha t}\Lambda_h({\tilde\Delta}_h^{-1}\hat\bta)| &\leq& C\bigg(\|\nabla\bu_h\|+\|\bu_h\|^{1/2}\|\nabla\bu_h\|^{1/2}
+\|\bu\|^{1/2}\|\nabla\bu\|^{1/2} \bigg)\bigg(\|\hat\bta\|_{-1}^{1/2}\hat\bta\|^{3/2}+\|\hat\bta\|\|\hat\bxi\| \bigg)\nonumber\\
&\leq& \epsilon\|\hat\bta\|^2+ C(\epsilon)\bigg(\|\nabla\bu_h\|^2+\|\bu_h\|\|\nabla\bu_h\|+\|\bu\|\|\nabla\bu\| \bigg)\|\hat\bxi\|^2
+ C(\epsilon)\|\hat\bta\|_{-1}^2\nonumber\\
&&\bigg(\|\nabla\bu_h\|^4+\|\bu_h\|^2\|\nabla\bu_h\|^2+\|\bu\|^2\|\nabla\bu\|^2 \bigg).
\een
Put $\epsilon = \frac{\nu}{2}$ in (\ref{L64.1***}) and use Lemmas \ref{L42} and \ref{L41} to obtain
\ben
\label{L64.2}
\frac{d}{dt}(\|\hat{\bta}\|_{-1}^2 +\kappa \| \hat{\bta}\|^2)+(\nu-\kappa\alpha)\|\hat{\bta}\|^2\leq 
C\|\hat\bxi\|^2+ (C(\nu)+2\alpha)\|\hat{\bta}\|_{-1}^2.
\een
Integrate (\ref{L64.2}) with respect to time and observe that $\bta(0)=0$ 
\ben
\label{L64.3}
\|\hat{\bta}\|_{-1}^2 +\kappa \| \hat{\bta}\|^2+(\nu-\kappa\alpha)\int_0^t\|\hat{\bta}\|^2ds\leq C\int_0^t\|\hat\bxi\|^2 ds
 + (C(K,\nu)+2\alpha)\int_0^t\|\hat{\bta}\|_{-1}^2 ds.  
\een
Apply Gronwall's Lemma in (\ref{L64.3}) and use Lemma \ref{L63}. Now, a use of triangular inequality completes the rest of proof.
\hfill {$\Box$}

\noindent
Now, we derive the main Theorem \ref{T41} of this section.
\begin{tdf}
\label{T41}
Let the assumptions {\rm (\bf {A1})-(\bf{A2})} 
and {\rm (\bf{B1})-(\bf{B2})} 
be satisfied. Further, let
the discrete initial velocity $\bu_{0h}= P_h \bu_0.$ 
Then, there exists a positive constant $C=C(\nu,~ \lambda_1,~\alpha, M)$ such that, 
for all $t\in (0,T]$ and for $\displaystyle{0\leq\alpha < \frac{\nu \lambda_1}{4\big(1+\lambda_1\kappa\big)}}$,
the following estimate holds:
\begin{equation}\label{velocity-estimate}
\|(\bu-\bu_h)(t)\|+ h\|\nab (\bu-\bu_h)(t)\|\le
C e^{C T} \kappa^{-1/2} h^2 .
\end{equation}
\end{tdf}
\noindent
{\it Proof.} Since $\e=\bu-\bu_h=(\bu-\bvh)+(\bvh-\bu_h)=\bxi+\bta$ and the estimate
of $\bxi$ is derived in Lemma \ref{L63}, therefore to complete the proof, it is enough to 
estimate $\bta$. 

With a choice of  ${\bphi}_h=e^{2\alpha t}\bta$ in (\ref{E653}), we apply  (\ref{2.1*}) to arrive at
\begin{equation}
\label{E654}
\frac{1}{2}\frac{d}{dt}(\|\hat{\bta}\|^2 +\kappa \|\nabla \hat{\bta}\|^2)+\big(\nu-\alpha (\kappa +\frac{1}{\lambda_1})\big) 
\|\nabla \hat{\bta}\|^2=e^{\alpha t}\Lambda_h(\hat{\bta}),
\end{equation}
where $\Lambda_h({\bphi}_h)$ is given as in (\ref{Lambda})). 
For the term on the right hand side of (\ref{E654}), we first rewrite it as
$$ 
e^{\alpha t}\Lambda_h(\hat{\bta})=e^{-\alpha t}\big(-b(\hat{\bf e}
,\hat{\bu}_h,\hat{\bta})+b(\hat{\bu},\hat{\bta},\hat{\bf e})\big).
$$
An application of the H\"older inequality with the Poincar\'{e} inequality, the Agmon inequality 
(\ref{sobolev-inequality}) 
and  the discrete  Sobolev inequality (see, Lemma 4.4 in \cite {HR82}) shows 
\begin{align}\label{Lambda-estimate}
e^{\alpha t}|\Lambda_h(\hat{\bta})|&\leq  C  e^{-\alpha t}\big(\|\hat{\bf
e}\|\|\nabla \hat{\bu}_h\|_{L^6}\|\hat\bta\|_{L^3}+\|\hat \bu\|_{L^\infty}\|\nabla\hat\bta\|\|\hat {\bf e}\|\big)\nonumber\\
&\leq  C\big( e^{-\alpha t} \|\tilde \Delta_h {\bu}_h \|\|\nabla\hat \bta\|\|\hat {\bf e}\|+\|\nabla \hat\bu\|^{\frac{1}{2}}\|\tilde{\Delta}\hat\bu\|^{\frac{1}{2}} \|\nabla \hat\bta\|\|\ \hat {\bf e}\|\big)\nonumber\\
&\leq  C\big(\nu) e^{-2\alpha t} (\|\tilde \Delta_h{\hat \bu_h}\|^2+\|\nabla \hat\bu\|\|\tilde \Delta\hat\bu\|\big)\|\hat {\bf e} \|^2+ \frac{\nu}{2}\|\nabla \hat\bta\|^2.
\end{align}
Substitute ${\bf e}=\bxi+\bta$ in  (\ref{Lambda-estimate}) to find that
\begin{eqnarray}\label{E656}
e^{\alpha t}|\Lambda_h(\hat{\bta})|\leq C(\nu) e^{-2\alpha t} (\|\tilde \Delta_h{\hat \bu_h}\|^2+ \|\nabla\hat\bu\|\|\tilde \Delta\hat\bu\|)( \|\hat\bxi\|^2+\|\hat\bta\|^2)+ \frac{\nu}{2}\|\nabla \hat\bta\|^2. 
\end{eqnarray}
A use of  (\ref{E656}) in (\ref{E654}) now yields
\begin{align}\label{E658*}
\frac{d}{dt}(\|\hat{\bta}\|^2 +\kappa \|\nabla \hat{\bta}\|^2)+(\beta +\nu) &\|\nabla \hat{\bta}\|^2  
\leq  C(\nu) e^{-2\alpha t}\big(( \|\hat\bxi\|^2+\|\hat\bta\|^2) \|\tilde \Delta_h{\hat \bu_h}\|^2\nonumber\\
  &+ ( \|\hat\bxi\|^2+\|\hat\bta\|^2)\|\nabla\hat\bu\|\|\tilde \Delta\hat\bu\|\big)
  + \nu \|\nabla \hat\bta\|^2.
\end{align}
Integrate (\ref{E658*}) with respect to time from $0$ to $t$ 
and  apply Lemmas \ref{L43-1}, \ref{L41} and \ref{L63} 
to arrive at 
\begin{align}
\|\hat{\bta}\|^2 +\kappa \|\nabla \hat {\bta}\|^2+ \beta &\int_0^t\|\nabla \hat{\bta}\|^2\;ds  \leq 
C(\nu,\alpha, \lambda_1, M) h^4 e^{2\alpha t}\nonumber\\
 &+\int_0^t \|\hat\bta\|^2(\|\nabla\bu\|\|\tilde \Delta\bu\|+\|\tilde \Delta_h{ \bu_h} \|^2)ds.
\end{align}
Then, use Gronwall's Lemma and then multiply by $e^{-2\alpha t}$ to obtain 
 \begin{eqnarray}\label{E670}
\|{\bta}\|^2+ \kappa\|\nabla \bta\|^2+ \beta e^{-2 \alpha t} \int_0^t{\|\nab\hat{\bta}(s)\|^2\,}ds\le
C h^4\exp{\left(\int_0^T(\|\nabla\bu\|\|\tilde \Delta \bu\|+
\|\tilde\Delta_h {\bu}_h\|^2)\,ds\right)}.
\end{eqnarray}
For the integral on the right hand side of (\ref{E670}), apply Lemmas \ref{L43} and \ref{L41} to arrive at
\begin{align}\label{E670*}
\int_0^T(\|\nabla\bu\|\|\tilde \Delta \bu\|+
\|\tilde\Delta_h {\bu}_h\|^2)\,ds\leq CT.
\end{align}
 Apply (\ref{E670*}) in (\ref{E670}) to derive estimates for $\bta$ as
\begin{eqnarray}\label{E671}
\|\bta\|^2+\kappa\|\nabla \bta\|^2+ 2\beta e^{-2\alpha t}\int_0^t e^{2\alpha s}\|\nab{\bta}(s)\|^2\,ds\le
C h^4 e^{CT}.
\end{eqnarray}
A use of triangle inequality along with (\ref{E671}) and Lemma \ref{L63} completes the rest of the 
proof.\hfill {$\Box$}
\begin{remark}
We observe that in the above proof the presence of the exponential 
term on the right-hand side of the error Theorem \ref{T41} is due to the 
estimate of $\bta$, as the estimate $ \bxi$ is uniform in time. 
In fact, the contribution of the exponential term comes from the Lemma \ref{L64}. If $\bu_0$ and ${\bf f}$ 
are sufficiently small with respect to the norms in the assumptions ({\bf A2}) so that
\ben \label{smallness}
\nu-(\kappa\alpha+C(K,\nu)+ 2\alpha) \ge 0.
\een   
then, from (\ref{L64.2}), we have
\ben
\label{nc8}
\frac{d}{dt}(\|\hat{\bta}\|_{-1}^2 +\kappa \| \hat{\bta}\|^2)+(\nu-(\kappa\alpha+(C(K,\nu)+2\alpha))\|\hat{\bta}\|^2\leq C(K,\nu)\|\hat\bxi\|^2 .\nonumber
\een
Integrate (\ref{nc8}) with respect to time $0$ to $t$ and use $\bta(0)=0$ to arrive at
\be
(\|\hat{\bta}\|_{-1}^2 +\kappa \|\hat{\bta}\|^2)+(\nu-(\kappa\alpha+(C(K,\nu)+2\alpha))\int_0^t\|\hat{\bta}\|^2 ds\leq  C(K,\nu)\int_0^t\|\hat\bxi\|^2 ds
\ee 
We can now avoid Gronwall's Lemma and  use Lemma \ref{L63} with triangle inequality to obtain
\be
e^{-2\alpha t}\int_0^t e^{2\alpha t}\|{\bf e}\|^2 ds \le C \kappa^{-1} \; h^4.
\ee
Following similar lines of proof, one can show the estimate of $\|\e(t)\|$ for all $t>0$ from Theorem \ref{T41}, 
provided the assumption (\ref{smallness}) is satisfied.
\end{remark}
\begin{remark}
 When $\f \in L^2(0,\infty; \bL^2(\Omega))$, all the error estimates are valid uniformly in time as
 all the {\it a priori} bounds hold true for $\alpha=0$ and therefore,
 the estimate (\ref{u-integral}) bounded uniformly in time. 
 Moreover,if $\f=0$ or $\f=O(e^{-\alpha_0 t}),$ we have as in \cite{ANS} exponential decay property
 for the solution as well as for the error estimates. 
\end{remark}

\noindent
{\bf Uniform in time estimates for the velocity}: We now derive uniform (in time) error estimate for 
the velocity term under the following uniqueness condition
\ben 
\label{uni} 
\frac{N}{\nu^2}\|{\bf f}\|_{L^{\infty}(0,\infty,\bL^2(\Omega)} < 1\;\;\;\; \text{and}\;\;\;\; N=\displaystyle\sup_{u,v,w \in V}\frac{|b(\bu,\bv,\bw)|}{\|\nabla \bu\|\|\nabla \bv\|\|\nabla \bw\|}. 
\een
When $f=0$ or $\|f(t)\|=\mathcal{O}(e^{-\alpha_0 t})$ for some $\alpha_0>0$, (\ref{uni}) satisfies trivially.
\begin{tdf}
\label{T51}
Under the assumption of Theorem \ref{T41} and the uniqueness condition (\ref{uni}), there exist a positive 
constant C, independent of time and $\kappa,$ such that for all $t >0$  
\ben \label{velocity-uniform}
\|(\bu-{\bu}_h)(t)\| + h\;\|\nabla(\bu-{\bu}_h)(t)\|\;\leq C\; \kappa^{-1/2}\; h^2. 
\een
\end{tdf}
\noindent
{\it Proof.} In order to derive estimates, which are valid uniformly for all $t >0$, 
we need derive a different estimate for the nonlinear term $ \Lambda_h(\hat{\bta})$ with the help of 
the uniqueness condition (\ref{uni}). Therefore, 
we rewrite
\ben
\label{T51.1}
\Lambda_h({\bta})=-[b(\bxi,\bu_h,\bta)+b(\bta,\bu_h,\bta)+b(\bu,\bxi,\bta)].
\een 
Using uniqueness condition, it follows that
\ben
\label{T51.2}
|b(\bta,\bu_h,\bta)| \leq N \|\nabla\bta\|^2\|\nabla\bu_h\|.
\een
Apply (\ref{L64.1*}) and (\ref{L64.1**}) to find that
\ben
\label{T51.3}
|b(\bxi,\bu_h,\bta)+b(\bu,\bxi,\bta)|\leq C\bigg(\|\tilde\Delta \bu\|^2+ \|\nabla \bu_h\|^{1/2}\|\tilde\Delta \bu_h\|^{1/2}\bigg)\|\nabla \bta\|\|\bxi\|.
\een
Substitute (\ref{T51.2}), (\ref{T51.3}) in (\ref{T51.2}) and use Lemma \ref{L63} to obtain
\ben
|\Lambda_h({\bta})| \leq N\|\nabla\bta\|^2\|\nabla\bu_h\|+Ch^2\|\nabla\bta\|.
\een  
Now, we modify the proof of Theorem \ref{T41} as follows
 \ben
\label{T51.4}
\frac{1}{2}\frac{d}{dt}(\|\hat{\bta}\|^2 +\kappa \|\nabla \hat{\bta}\|^2)+(\nu-N\|\nabla \bu_h\|)\|\nabla\hat{\bta}\|^2
\leq \alpha (\|\hat{\bta}\|^2+\kappa \|\nabla \hat{\bta}\|^2)+ Ch^2 \|\nabla\hat{\bta}\|.
\een
An integration with respect to time with multiplication by $e^{2\alpha t}$ leads to 
\begin{align}
\label{T51.5}
\|{\bta}(t)\|^2 &+\kappa \|\nabla{\bta}(t)\|^2+2e^{-2\alpha t}\int_0^t e^{2 \alpha s}(\nu-N\|\nabla \bu_h\|)
\|\nabla {\bta}(s)\|^2 \nonumber\\
&\leq 2\alpha e^{-2\alpha t}\int_0^t e^{2 \alpha s} (\|{\bta}(s)\|^2
+\kappa \|\nabla {\bta}(s)\|^2) ds+ Ch^2 e^{-2\alpha t}\int_0^t e^{2 \alpha s}\|\nabla{\bta}(s)\| ds.
\end{align}
Letting $t \rightarrow \infty$, we obtain
\ben
\frac{1}{\nu}\bigg(1-N \nu^{-2}\|{\bf f}\|_{\bL^{\infty}(0,\infty,\bL^2(\Omega)}\bigg) 
\limsup_{t \rightarrow \infty}\|\nabla\bta(t)\| \leq C h^2.
\een
Then, we conclude from the uniqueness condition (\ref{uni}) that
\ben 
\limsup_{t \rightarrow \infty}\|\nabla\bta(t)\| \leq C h^2,
\een
and hence,
\ben \label{eta-sup}
\limsup_{t \rightarrow \infty}\|\bta(t)\| \leq C h^2.
\een
Now the uniform estimate of $\bxi$  combined with (\ref{eta-sup}) leads to 
\ben
\limsup_{t \rightarrow \infty}\|{\bf e}(t)\| \leq C\;\kappa^{-1/2}\; h^2.
\een
Note that C is valid uniformly for all $t > 0$, and this complete the rest of the proof.
\hfill {$\Box$} 
\section{ \normalsize \bf Error estimate for the pressure}
\setcounter{tdf}{0}
\setcounter{ldf}{0}
\setcounter{cdf}{0}
\setcounter{equation}{0}
In this section, the optimal error estimate for the
Galerkin approximation $p_h$ of the pressure $p$ is derived. Further, under the uniqueness condition (\ref{uni}), 
the estimate is shown to be valid uniformly in time. The main theorem of this section is stated as follows: 
\begin{tdf}
\label{Th61}
Under the hypotheses of Theorem \ref{T41},  
there exists a positive constant $C$ depending on $\nu,~ \lambda_1,~\alpha$ and $ M$, 
such that for $T>0$ with $0<t\leq T$
$$
\|(p-p_h)(t)\|_{L^2/N_h}\leq C e^{C T} \kappa^{-1/2}\; h .
$$
\end{tdf}  
\noindent
We prove the  theorem \ref{Th61} with help of Lemmas \ref{L71} and \ref{L73}. From ({\bf{B2}}), it follows
that
\ben
\label{T61.1}
\|(j_h p-p_h)(t)\|_{L^2/N_h}\leq C\bigg(\|j_h p-p\|+\displaystyle{\sup_{{\bphi}_h\in\bH_h/\{0\}}}{\left\{\frac{( p-p_h,\nab\cdot{\bphi}_h)}{\|\nab{\bphi}_h\|}\right\}}\bigg).
\een
We observe that  the estimate of the first term on the right hand side of (\ref{T61.1})
follows from the approximation property stated in ({\bf B1}). To
complete the proof, it is sufficient to estimate the
second term in (\ref{T61.1}). 
Use (\ref{4.1}) and  (\ref{E62}) to find that for $\bphi_h \in \bH_h$
$$
(p-p_h,\nab\cdot{\bphi}_h)=({\bf e}_t,{\bphi}_h)+\kappa\,a({\e_t},{\bphi}_h)+ \nu a({\e},{\bphi}_h)-
\Lambda_h({\bphi}_h)\,\,\, \forall {\bphi}_h\in\bH_h,
$$
where $ \Lambda_h({\bphi}_h)$ is given as in (\ref{Lambda}).
A use of generalized H\"olders inequality  with Sobolev imbedding, Lemmas \ref{L42} and \ref{L41} leads to
\begin{eqnarray}\label{E71*}
|\Lambda_h({\bphi}_h)|\le C(\|\nabla\bu_h\|+\|\nabla\bu\|)\|\nab{\bf
e}\|\|\nab{\bphi}_h\| \leq C \|\nab{\bf e}\|\|\nab{\bphi}_h\|. 
\end{eqnarray}
Thus,
$$
(p-p_h,\nab\cdot{\bphi}_h)\le C(\nu)\big(\|{\bf e}_t\|_{-1,h}
+\kappa \|\nab{\bf e}_t\|+\kappa \|\nab{\bf e}\|  \big)\|\nab{\bphi}_h\|,
$$
where 
$$\|{\bf e}_t\|_{-1,h}= \displaystyle{\sup_{{\bphi}_h\in\bH_h/\{0\}}} 
{\left\{\frac{({\bf e}_t, {\bphi}_h)}{\|\nab{\bphi}_h\|}\right\}}.
$$
Altogether, we derive the following result.
\begin{ldf} \label{L71}
 The semidiscrete Galerkin approximation $p_h$ of
the pressure $p$ satisfies for all $t \in (0, T]$ 
\begin{eqnarray}
\label{E72}
\|(p-p_h)(t)\|_{L^2/N_h}\le C\big(\|{\bf e}_t\|_{-1,h}+ \kappa \|\nab{\bf e}_t\|+ \|\nab{\bf e}\|\big). 
\end{eqnarray}
\end{ldf}
\noindent
Note that  the estimate $\|\nab{\bf e}\|$ is known from the  Theorem \ref{T61}. In order to complete 
the proof of Theorem \ref{Th61}, we only need to estimate $\|{\bf e}_t\|_{-1,h}$ and $\|\nabla {\bf e}_t\|$.
\begin{ldf}\label{L73}
For all $t\in (0,T]$, the error ${\bf e}=\bu-\bu_h$ in the velocity satisfies 
\begin{eqnarray}
\|{\bf e}_t(t)\|_{-1,h}+\kappa \|\nabla {\bf e}_t(t)\|
 \leq C e^{CT}\kappa^{-1/2}\;h.
\end{eqnarray}
\end{ldf}
\noindent{\it Proof.} Subtract (\ref{4.2}) from (\ref{E62}) to write
\begin{eqnarray}
\label{E74*}
({\bf e}_t,{\bphi}_h)+\kappa\,a({\bf e}_t, {\bphi}_h)+\nu a({\bf e}, {\bphi}_h)=\Lambda_h({\bphi}_h)+(p,\nab\cdot{\bphi}_h),\,\,{\bphi}_h\in \bH_h.
\end{eqnarray}
 where $\Lambda_{h}({\bphi}_h)$ is defined in (\ref{Lambda}).
Choose ${\bphi}_h=P_h{\bf e}_t= \e_t-(\e_t-P_h \e_t)=  \e_t-(\bu_t-P_h\bu_t)$ in (\ref{E74*}) to arrive at
\begin{eqnarray}
\label{E72*}
 \|{\bf e}_t\|^2+\kappa\,\|\nabla{\bf e}_t\|^2 &=&({\bf e}_t, \bu_t-P_h \bu_t)+\kappa\,a({\bf e}_t, \bu_t-P_h \bu_t)+(p, \nabla\cdot  P_h{\bf e}_t)\nonumber\\
 &&+ \Lambda_{h}(P_h{\bf e}_t)-\nu  a({\bf e},P_h \e_t).
\end{eqnarray}
For the first term together with the last term on the right hand side of (\ref{E74*}), 
apply Poincar\'e inequality (\ref{2.2}) and the stability property of $P_h$ to 
obtain
\ben \label{E77-1}
({\bf e}_t, \bu_t-P_h \bu_t) -\nu  a({\bf e},P_h \e_t)  \leq \Big(\lambda_1^{-1/2} \; \|\bu_t-P_h \bu_t\|
+ \nu \|\nabla \e\|\Big)\;\|\nabla \e_t\|.
\een
For the third term on the right hand side of (\ref{E74*}), a use of
the discrete incompressible condition with (\ref{4.4}) yields
\ben
\label{E77*}
|(p,\nab\cdot P_h{\bf e}_t)|=|(p-j_h p,\nab\cdot P_h{\bf e}_t)|\le
\|p-j_h p\|\|\nab{\bf e}_t\|.
\een
In order to estimate the fourth term on the right hand side of (\ref{E74*}),  
apply (\ref{E71*}) and (\ref{4.4}) to obtain 
\begin{eqnarray}
\label{E76*}
|\Lambda_{h}(P_h{\bf e}_t)|\leq C \|\nabla {\bf e}\|\|\nabla {\bf e}_t\|.
\end{eqnarray}
Substitute (\ref{E77-1}), (\ref{E77*}) and (\ref{E76*}) in (\ref{E72*}) to arrive at
\begin{eqnarray}
\label{E78*}
\kappa\|\nabla {\bf e}_t\| \leq  C(\nu,\lambda)\big( \|\nabla {\bf e} \|
+ \kappa \|\nabla( \bu_t-P_h \bu_t)\|+\|p-j_h p\| +\|\bu_t-P_h \bu_t\| \big).
\end{eqnarray}
A use of  (\ref{4.5}) and ({\bf B1}) in (\ref{E78*}) shows
\begin{equation}
\label{E79*}
\kappa\|\nabla {\bf e}_t\| \leq C(\nu,\lambda)\big(\|\nabla{\bf e}\|
+h ( \kappa \|\tilde\Delta\bu_t\| +\|\nabla p\| + \kappa^{-1/2}\; \kappa^{1/2}\|\nabla \bu_t\| )\big).
\end{equation}
An application of Theorems \ref{T31} and \ref{T61} with  \ref{estimate-Delta-1}  shows that
\begin{equation}
\label{E79-1}
\kappa\|\nabla {\bf e}_t\| \leq C(\nu,\lambda,\alpha, M) \kappa^{-1/2}\; h.
\end{equation}
To complete the rest of the proof, observe from (\ref{E63})that
\begin{eqnarray}\label{E63-p}
({\bf e}_{t},{\bphi}_h)=-\kappa a({\bf e}_{t},{\bphi}_h)-\nu a({\bf e},{\bphi}_h)+ {\bf \Lambda}({\bphi}_h)
+(p,\nabla \cdot {\bphi}_h)
\end{eqnarray}
An application of  the Cauchy-Schwarz inequality to (\ref{E63-p}) with estimates (\ref{E77*}) and (\ref{E76*}) 
shows
\begin{eqnarray}\label{E63-p-1}
({\bf e}_{t},{\bphi}_h) \leq \Big( \kappa \|\nabla {\bf e}_{t}\| +\nu \|\nabla{\bf e}\|+ C \|\nabla \e\| 
+ \|(p-j_h p)\|\Big)\;\|\nabla {\bphi}_h\|,
\end{eqnarray}
and hence, a use of ({\bf B1}) with theorem \ref{T61} and estimate (\ref{E79-1}) yields
the estimate of $\|\e_t\|_{1,h}.$ This 
concludes the proof.\hfill {$\Box$}\\

\noindent {\it Proof of Theorem \ref{Th61}}.
The proof follows from Lemmas \ref{L71} and \ref{L73} with the approximation property $({\bf B1})$ 
of $j_h$.\hfill {$\Box$}

\begin{remark}
Under uniqueness condition (\ref{uni}), an appeal to (\ref{E72}) and (\ref{E79*}) leads to
the error estimate for the pressure,  which is valid for all time $t >0$: 
\ben \label{p-uniform}
\|(p-p_h)(t)\|_{L^2/N_h}\leq K\; \kappa^{-1/2}\; h, 
\een
and this provides optimal error estimate for pressure term, which is valid uniformly in time.
\end{remark}

\begin{remark}
In Theorems \ref{T41}, \ref{T51} and \ref{T61}, if we choose $\kappa^{1/2} =O(h^{\delta}),$ where $\delta>0$ 
can be take sufficiently small, then we obtain the following quasi-optimal order of convergence:
\ben \label{error}
\|(\bu-{\bu}_h)(t)\| + h\;\Big(\|\nabla(\bu-{\bu}_h)(t)\| + \|(p-p_h)(t)\|_{L^2/N_h}\Big)= O(h^{2-\delta}). 
\een
\end{remark} 

\section{\bf{Numerical Experiments}} 
 \setcounter{tdf}{0}
\setcounter{ldf}{0}
\setcounter{cdf}{0}
\setcounter{equation}{0} 
In this section, three numerical examples using mixed finite element space $P_2$-$P_0$ for 
spatial discretization and backward Euler method 
for temporal discretization are discussed with computed orders of convergence, which confirm  our theoretical 
findings. Moreover, it is shown through numerical experiments that orders of convergence do not deteriorate with 
$\kappa$ small which again matches with theory.
For all three examples, consider the domain $\Omega=(0,1)\times(0,1),$ $T=1,$ $\kappa=$ and $\nu=1.$ 
Choose approximating spaces 
$\bH_h$ and $L_h$ for velocity and pressure, respectively,  as  
\begin{eqnarray}
{\bf H}_h= \{{\bf v}\in {\big( C(\bar\Omega)}\big)^2: {\bf v}|_K \in (P_2(K))^{2},K \in \tau_h \}\;\;
\mbox{and }\;
L_h= \{q \in L^2(\Omega):  q|_K \in P_0(K),K \in \tau_h \},\nonumber
\end{eqnarray}
where $\tau_h$ denotes an admissible triangulation of $\bar\Omega$ in to closed triangles with mesh size $h.$
Let $0=t_0<t_1< \cdots< t_N=T$, be a uniform subdivision of the time interval $(0, T]$ with 
$t_{n}=nk$  and $k=t_n-t_{n-1}.$  
The fully discrete backward Euler method can be formulated as: 
given ${\bU}^{n-1}$, find the pair $(\bU^n,P^n)$ approximating the pair $({\bf u}, p)$ at $t=t_{n}=nk $ 
satisfying
\begin{eqnarray}\label{E84}
({\bar\partial_t}{\bf U}^{n},{\bf v}_h)+ \kappa a( \bar\partial_t {\bf U}^{n},{\bf v}_h)
&+& \nu a( \bU^n, {\bf v}_h)  +  b( {\bf U}^{n}, {\bf U}^{n},{\bf v}_h)
+ ( {\bf v}_h,\nabla P^n)\\
&=&({\bf f}^n,{\bf v}_h),\,\,\,\forall {\bf v}_h \in {\bf H}_h,\nonumber\\
(\nabla \cdot {\bf U}^n ,w_h) & = & 0 ,\,\,\,\,\,\,\,\,\,\,\, \forall  w_h \in  L_h,\nonumber
\end{eqnarray}
where $\bar\partial_t \bU^n=\displaystyle{\frac{\bU^n-\bU^{n-1}}{k}}$.
\noindent
\begin{example}\label{ex1}

The convergence rates of the approximate solution is verified by choosing the right hand side function $f$ 
in such a way that the exact solution $(\bu,p)=((u_1,u_2),p)$ of (\ref{1.1})-(\ref{1.3}) is  given as
$$
u_1  =   10 \cos t\,\, x^2(x-1)^2 y (y-1) (2y-1),\;\;\;\;u_2  =  -10 \cos t\,\, y^2(y-1)^2 x (x-1) (2x-1),\;\;\;\;p = 40\cos t\,\, xy.
$$
\end{example}
\noindent	
The theoretical analysis proves the convergence rates $\mathcal{O}(h^2)$ for velocity in $\bL^2$ norm, $\mathcal{O}(h)$ for velocity 
in $\bH^1$-norm and $\mathcal{O}(h)$ for pressure in $L^2$ norm. 
Figure 1 provides convergence rates obtained on successively 
refined meshes with time step size $k= h^2.$ These results agree with the optimal theoretical 
convergence rates obtained in Theorems \ref{T41}-\ref{Th61}.
\begin{figure}[ht!]
 \begin{minipage}[b]{.5\linewidth} 
\centering
{\includegraphics[width=2.5in]{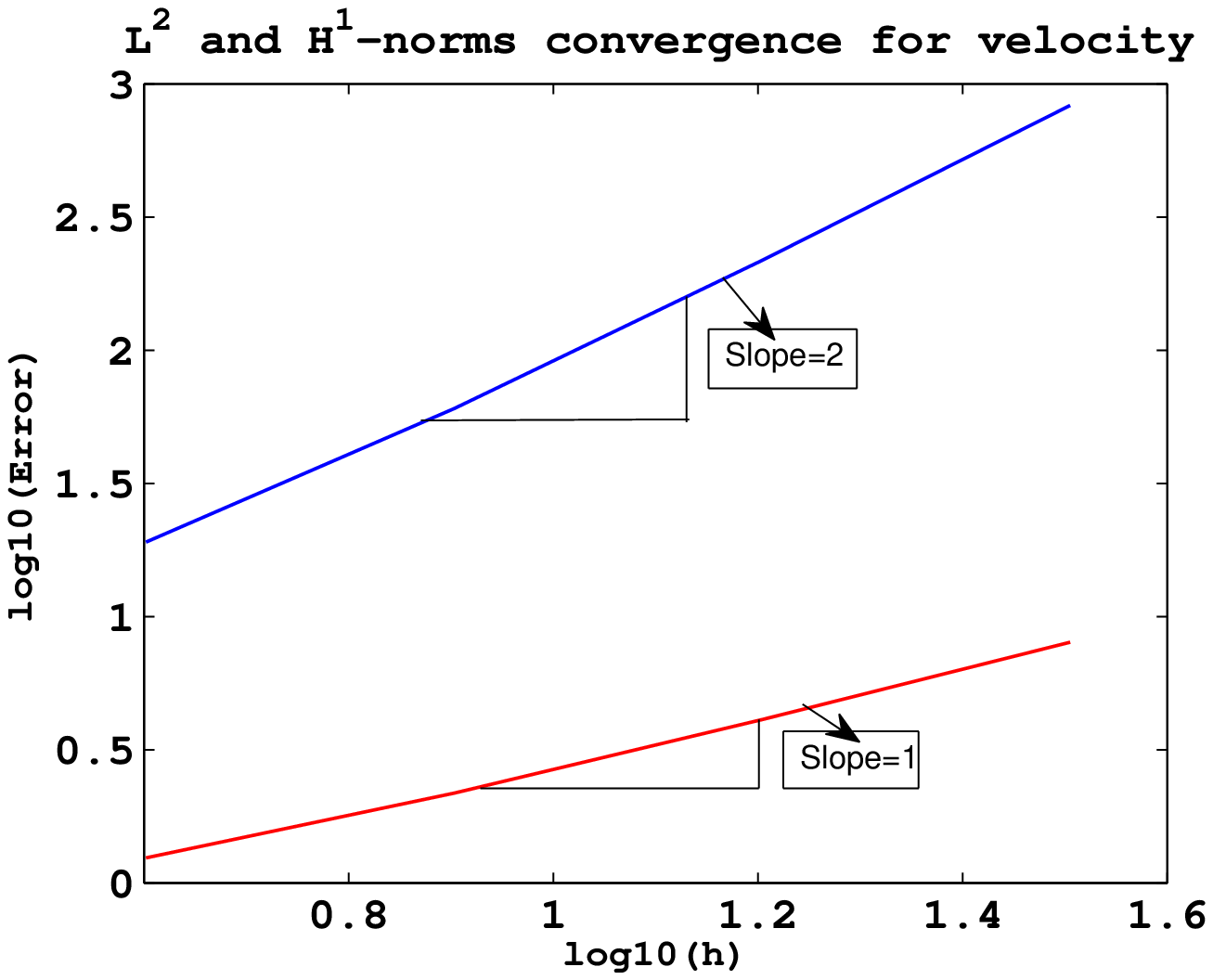}}
\end{minipage}
 \hspace{0.05cm} 
   \begin{minipage}[b]{0.5\linewidth}
 \centering
  { \includegraphics[width=2.5in]{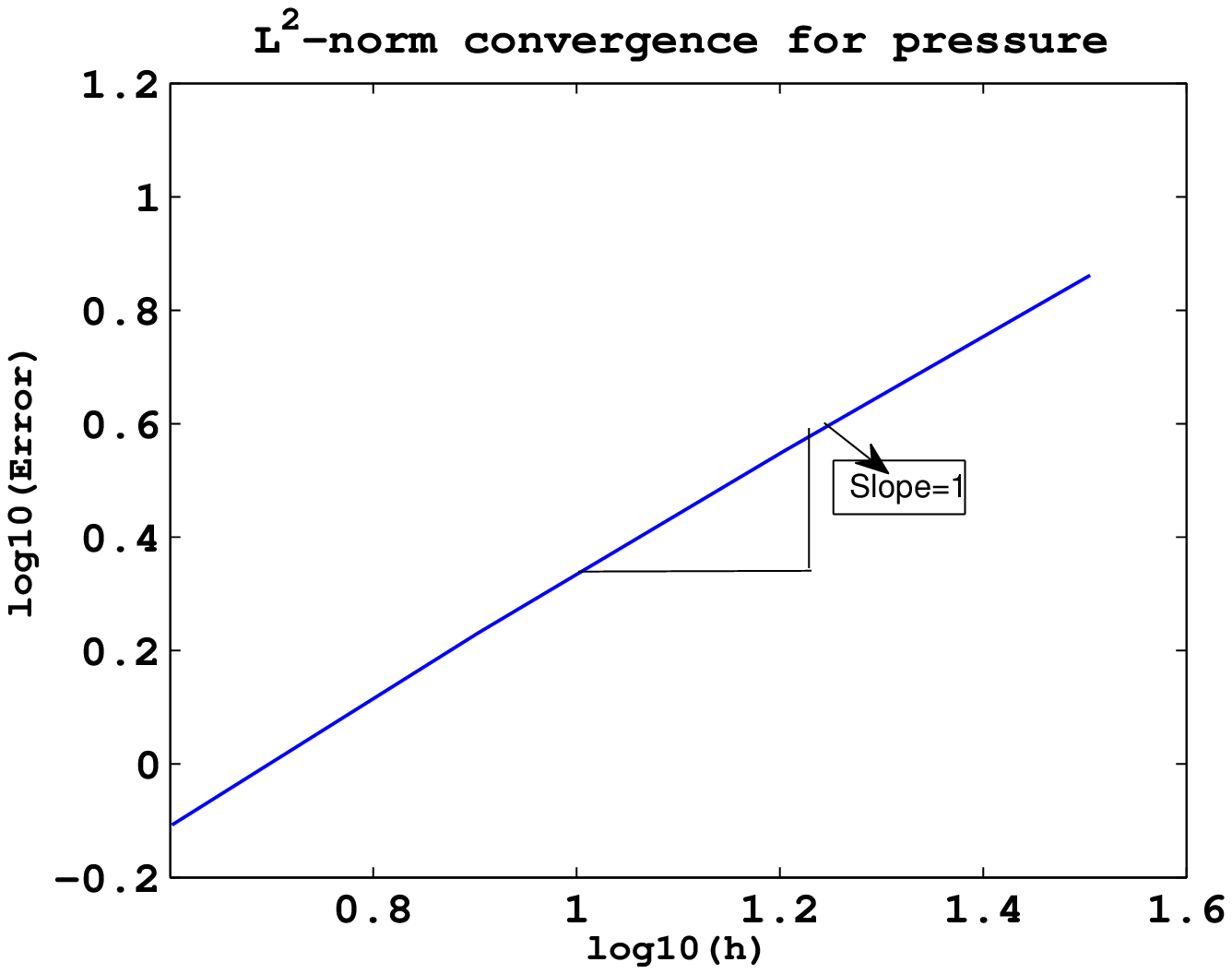}}
\end{minipage}
   \caption{Plot of Convergence Rates for Example \ref{ex1}}
   \end{figure}
    \begin{figure}[ht!]
 \begin{center}
 {\includegraphics[width=3.0in]{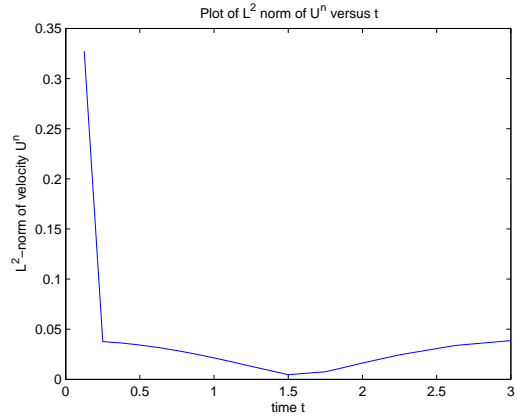}}
 \end{center}
 \caption{Boundedness of $\|\bU^n\|$ as time varies.}
    \end{figure}
\noindent
Figure 2 depicts that the approximate solution for the data in Example \ref{ex1} is bounded. Note that, here the right hand side 
function is bounded for all time. Further, Tables 1, 2, 3 represent that the order of convergence for the velocity and 
pressure errors in Theorems 5.1 and 6.1 hold true in the limit $\kappa \rightarrow 0$. 
\begin{table}[!ht]
\centering
\begin{tabular}{|c|c|c|c|c|c|c|c|}
\hline
S No    & $h$ & $\|\bu(t_n)-\bU^n\|_{\bL^2}$ & $\|\bu(t_n)-\bU^n\|_{\bL^2}$  & $\|\bu(t_n)-\bU^n\|_{\bL^2}$ 
& $\|\bu(t_n)-\bU^n\|_{\bL^2}$ \\
        &     & $\kappa=1$        &   $\kappa=10^{-3}$       & $ \kappa=10^{-6}$& $\kappa=10^{-9}$ \\
\hline
\hline1 &  1/4   &  1.28476 & 1.46678& 1.46699 & 1.46699 \\
\hline2 &  1/8   &  1.66634 & 1.71546& 1.71552 & 1.71552 \\
\hline3 &  1/16  &  1.84754 & 1.86060& 1.86062 & 1.86062 \\
\hline4 &  1/32  &  1.93052 & 1.93390& 1.93391 & 1.93391 \\
\hline
\end{tabular}
\caption{Numerical convergence rates for velocity error with variation in $\kappa$ for Example \ref{ex1}}
\end{table}

\vspace{2cm}
\begin{table}[!ht]
\centering
\begin{tabular}{|c|c|c|c|c|c|c|c|}
\hline
S No    & $h$ & $\|\bu(t_n)-\bU^n\|_{\bH^1}$ & $\|\bu(t_n)-\bU^n\|_{\bH^1}$  & $\|\bu(t_n)-\bU^n\|_{\bH^1}$ 
& $\|\bu(t_n)-\bU^n\|_{\bH^1}$ \\
        &     & $\kappa=1$        &   $\kappa=10^{-3}$       & $ \kappa=10^{-6}$       & $\kappa=10^{-9}$\\
\hline
\hline1 &  1/4   &  0.52668 & 0.70916& 0.70938 & 0.70938\\
\hline2 &  1/8   &  0.80620 & 0.85510& 0.85516 & 0.85516\\
\hline3 &  1/16  &  0.91745 & 0.93032& 0.93033 & 0.93033\\
\hline4 &  1/32  &  0.96385 & 0.96716& 0.96716 & 0.96716\\
\hline
\end{tabular}
\caption{Numerical convergence rates for velocity error with variation in $\kappa$ for Example \ref{ex1}}
\end{table}
\vspace{2cm}
\begin{table}[!ht]
\centering
\begin{tabular}{|c|c|c|c|c|c|c|c|}
\hline
S No    & $h$ & $\|p(t_n)-P^n\|$ & $\|p(t_n)-P^n\|$  & $\|p(t_n)-P^n\|$ & $\|p(t_n)-P^n\|$ \\
        &     & $\kappa=1$        &   $\kappa=10^{-3}$       & $ \kappa=10^{-6}$       & $\kappa=10^{-9}$      \\
\hline
\hline1 &  1/4   & 1.25307 & 1.25165& 1.25164 & 1.25164 \\
\hline2 &  1/8   & 1.12462 & 1.11394& 1.11393 & 1.11393 \\
\hline3 &  1/16  & 1.06496 & 1.05938& 1.05937 & 1.05937 \\
\hline4 &  1/32  & 1.02882 & 1.02663& 1.02663 & 1.02663\\
\hline
\end{tabular}
\caption{Numerical convergence rates for pressure error with variation in $\kappa$ for Example \ref{ex1}}
\end{table}
\begin{example}\label{ex2}
In this example, the initial velocity is chosen as\\
$u_1  =   10 \,\, x^2(x-1)^2 y (y-1) (2y-1),\;\;\;\;u_2  =  -10 \,\, y^2(y-1)^2 x (x-1) (2x-1),\;\;\;\;p = 40\,\, xy$\\
\noindent
 with $\nu=1$, $\kappa=1$ and $f=0$. In this case, to obtain the error estimates the exact solution 
 $\bu$ is replaced  by finite element solution obtained in a refined mesh. 
\end{example}
\noindent
The convergence rates presented in Figure 3 are in agreement with the 
results obtained for $f=0$, that is, the convergence rate for velocity in $\bL^2$ norm is $\mathcal{O}(h^2)$, 
for velocity 
in $\bH^1$-norm is $\mathcal{O}(h)$ and for pressure in $L^2$ norm is $\mathcal{O}(h)$. In Figure 4, the exponential decay property for the approximate solution $\|\bU^n\|$ is 
shown which verifies theoretical estimates for $f=0$.
\begin{figure}[ht!]
 \begin{minipage}[b]{.5\linewidth} 
\centering
{\includegraphics[width=2.5in]{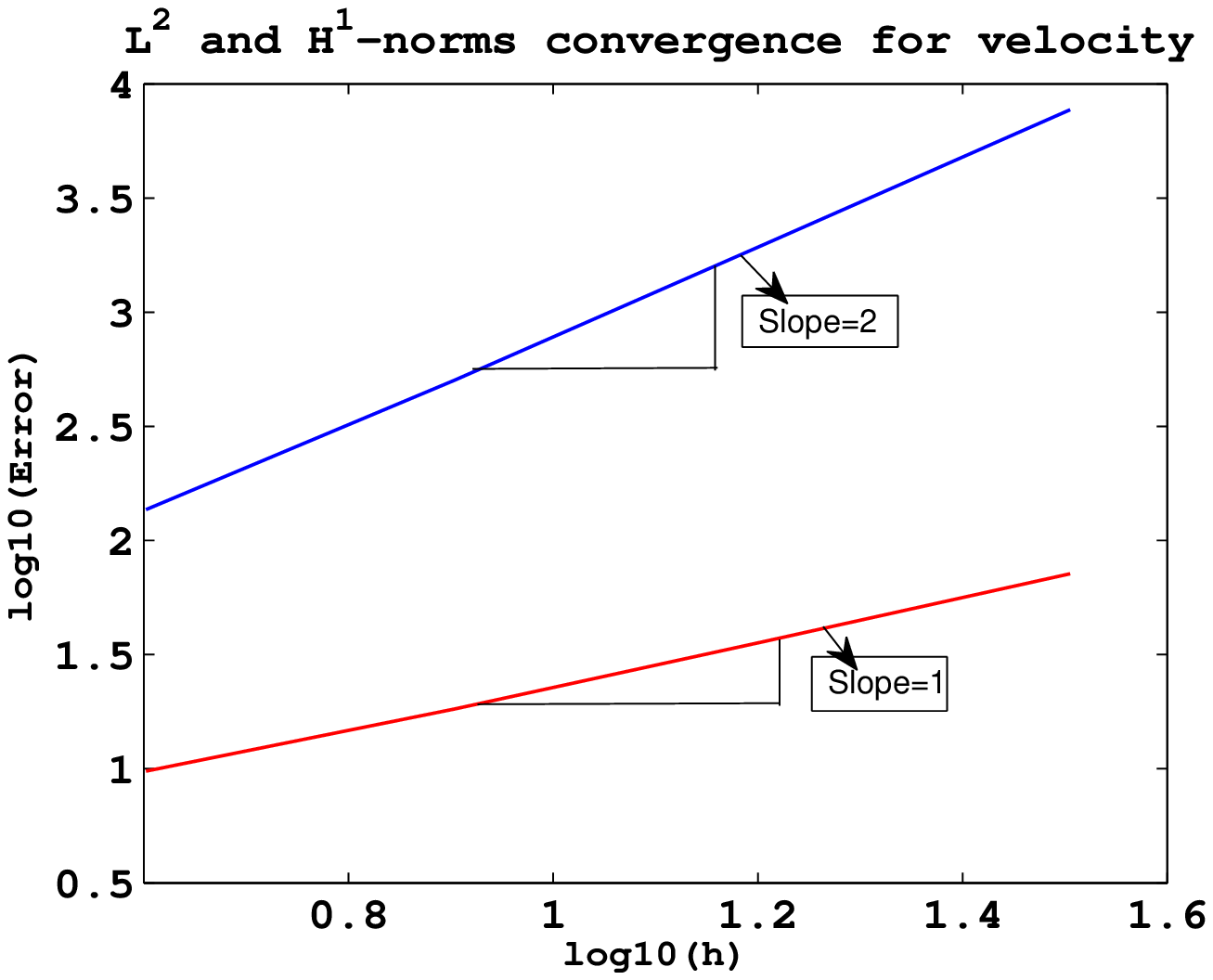}}
\end{minipage}
 \hspace{0.05cm} 
   \begin{minipage}[b]{0.5\linewidth}
 \centering
  { \includegraphics[width=2.5in]{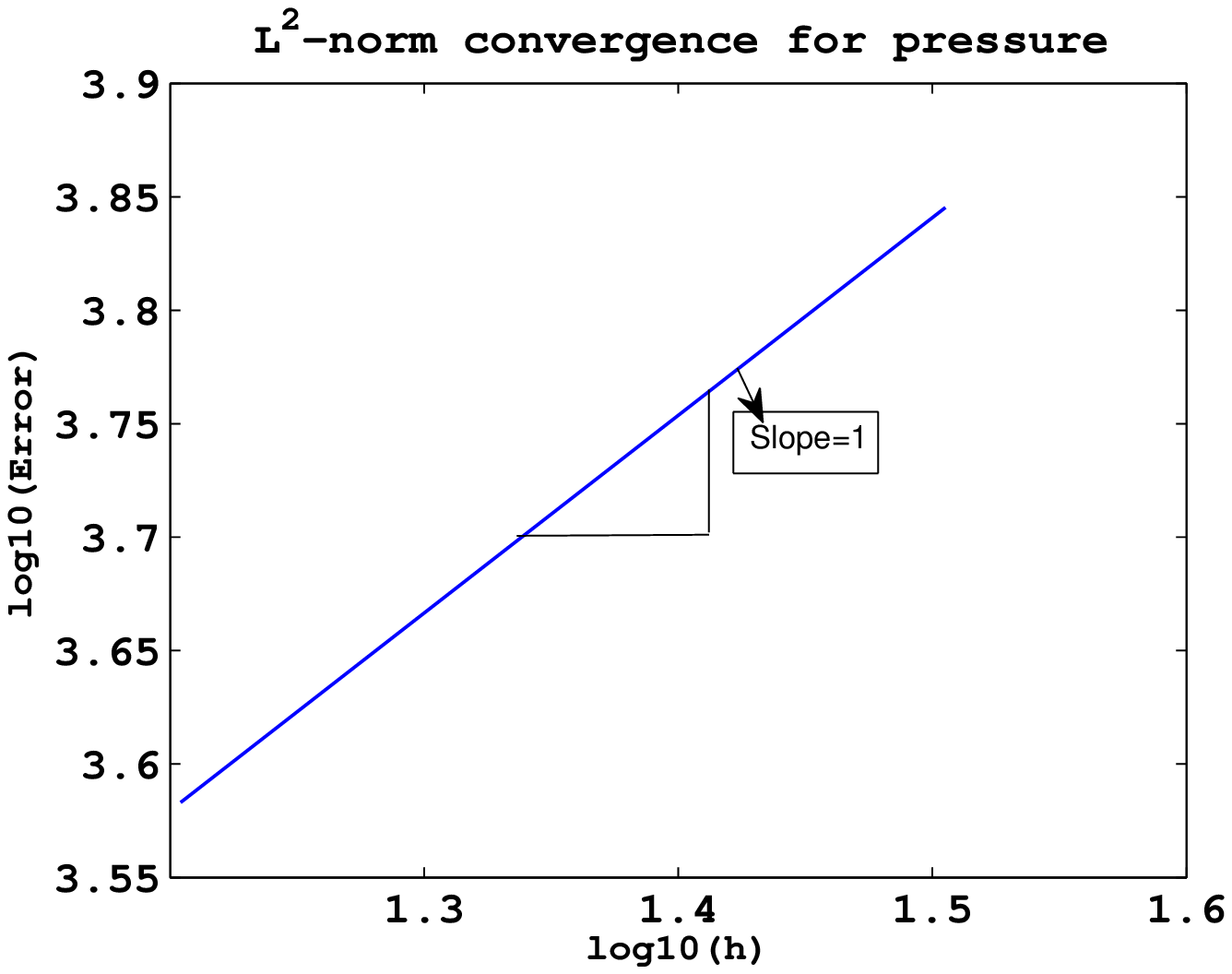}}
\end{minipage}
   \caption{Plot of Convergence Rates for Example \ref{ex2}}
   \end{figure}
   \vspace{2cm}
 \begin{figure}[ht!]
 \label{fg4}
 \begin{center}
 {\includegraphics[width=2.5in]{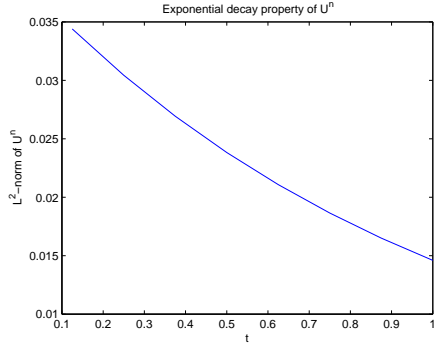}}
 \end{center}
 \caption{Decay property of $\|\bU^n\|$ corresponding to Example \ref{ex2}}
\end{figure}
\noindent
\begin{example}\label{ex3}
 This example demonstrates the exponential decay property of the discrete solution. Here, 
 $\nu=1$, $\kappa=1$ and $f=0$ with $\bu_0=\left(\sin^2(3\pi x)\sin(6\pi y), -\sin^2(3\pi y)\sin(6\pi x),sin(2\pi x)\sin(2\pi y)\right)$ in 
(\ref{1.1})-(\ref{1.3}). Once again, the error estimates are achieved by considering refined finite element solution as an exact solution. 
\end{example}
\noindent
The order of convergence is shown in Table 4. 
Figure 5 represents the exponential decay property of $\|\bU^n\|$ as time varies which is expected from theoretical analysis 
for right hand side function $f=0$.\\
\begin{table}[!ht]
\centering
\begin{tabular}{|l|l|l|l|l|l|l|l|}
\hline
S No & $h$ & $\|\bu-\bU^n\|_{\bL^2}$ & Convergence  & $\|\bu-\bU^n\|_{\bH^1}$ &  Convergence  \\
     &     &                         &    Rate      &                         &   Rate         \\
\hline
\hline1 &  1/4   & 0.430939 &                   & 6.833152943841204&         \\
\hline2 &  1/8   & 0.203398 & 1.083175531775576 & 5.967502741440636& 0.195424 \\
\hline3 &  1/16  & 0.065544 & 1.633758732735566 & 3.674410879988224& 0.699614\\
\hline4 &  1/32  & 0.017502 & 1.904904362752530 & 1.917811790943292& 0.938051\\
\hline
\end{tabular}
\caption{Numerical errors and Convergence rates with $k=h^2$ for Example \ref{ex3}}
\end{table}
\begin{figure}[ht!]
\label{fg5}
\centering
{\includegraphics[width=2.5in]{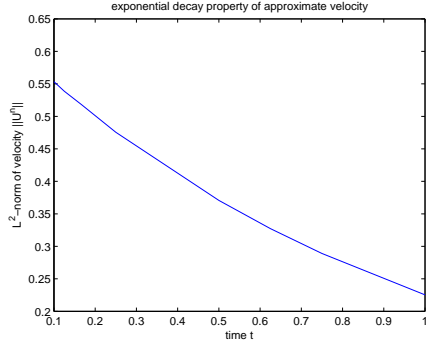}}
\caption{Decay property of $\|\bU^n\|$ corresponding to Example \ref{ex3}}
\end{figure}

\section{\bf{Conclusion}} 
This article in its first part deals with {\it a priori} estimates for the weak solution of 
(\ref{1.1})-(\ref{1.3}) which are valid uniformly in time  as  $t \mapsto \infty$ and also
uniformly for all $\kappa$ as $\kappa \mapsto 0.$ While estimates hold for $2D,$ that is, $d=2,$ and for $3D,$
that is, $d=3,$ estimates are valid  with smallness assumption on the data.In the second part, 
semidiscrete optimal error estimates  of order 
$O(\kappa^{-1/2} h^m)$ are derived for the velocity in $L^{\infty}(\bL^2)$-norm when $m=2$ and 
for the velocity in $L^{\infty}(\bH^1_0)$-norm, when $m=1$. Moreover for the pressure term, optimal order
estimate $L^{\infty}(L^2)$-norm, which is of order $O(\kappa^{-1/2} h)$ is established. In all these error
analyses, constants appeared in the error estimates depend exponentially on $T.$ But, under the uniqueness
assumption, it is shown that optimal error estimates are valid uniformly for all time $t>0.$ Further,
with $\kappa=O(h^{2\delta}),\;\delta>0$ very small, quasi-optimal error estimates are derived which are valid
uniformly in $\kappa$ as $\kappa\mapsto 0.$ All the above results  hold true for $2D$, but for $3D$ with smallness
assumption on the data. However, in stead of applying Lemma \ref{L43-1}, if  we apply 
Lemma \ref{L43}, then regularity results like in Theorem \ref{T31} can be obtained 
now with constants depending on $1/\kappa,$ but all results are valid for $3D$ without assumption of smallness
on the data. Similar conclusion for optimal error estimates can be derived, but with constants depending 
on $1/\kappa.$ Finally, numerical experiments are conducted to confirm our theoretical 
findings. 

\vspace{2em}
\noindent
{\bf{Acknowledgement}} The authors thank the referee for his/her valuable suggestion which has helped to 
improve our results.


\end{document}